\documentclass[11pt]{amsart}

\usepackage{eucal}

\usepackage{fullpage}

\def\hok{\mbox{}\begin{picture}(10,10)\put(1,0){\line(1,0){7}}
  \put(8,0){\line(0,1){7}}\end{picture}\mbox{}}

\newcommand{\lap}{\Delta}
\def\<{\langle}
\def\n{\nabla}
\def\>{\rangle}
\def\a{\alpha}

\def\j{\jmath}
\def\be{\beta}
\def\b{\bar}
\def\w{\wedge}
\def\o{\omega}
\def\O{\Omega}
\def\v{\varphi}
\def\s{{\sigma}}
\def\d{\partial}
\def\db{\overline{\partial}}
\def\ddb{\partial  \overline{\partial}}

\newtheorem{theorem}{Theorem}

\newtheorem{lemma}[theorem]{Lemma}
\newtheorem{proposition}[theorem]{Proposition}
\newtheorem{corollary}[theorem]{Corollary}
\newenvironment{remark}{\medskip \refstepcounter{theorem}
\noindent {\bf Remark \thetheorem .}}{\qed}
\newenvironment{conjecture}{\medskip \refstepcounter{theorem}
\noindent {\bf Conjecture \thetheorem .}}

\begin{document}
\title{Heat Flows for Extremal K\"ahler Metrics}
\author{Santiago R. Simanca}
\address{Institute for Mathematical Sciences, Stony Brook, NY 11794}
\email{santiago@math.sunysb.edu}

\begin{abstract}
Let $(M,J,\Omega)$ be a closed polarized complex manifold of K\"ahler type. 
Let $G$ be the maximal compact subgroup of the automorphism group of
$(M,J)$. On the space of K\"ahler metrics that are invariant under 
$G$ and represent the cohomology class $\Omega$, we define a flow
equation whose critical points are extremal metrics, those that  
minimize the square of the $L^2$-norm of the scalar curvature. We prove 
that the dynamical system in this space of metrics defined by the said 
flow does not have periodic orbits, and that its only fixed points
are extremal metrics. We prove local time existence of
the flow, and conclude that if the lifespan of the solution is finite, then
the supremum of the norm of its curvature tensor must blow-up as time
approaches it. While doing this, we also prove that extremal solitons can only
exist in the non-compact case, and that the range of the holomorphy potential 
of the scalar curvature is an interval independent of the metric chosen to
represent $\O$. 
We end up with some conjectures concerning the plausible
existence and convergence of global solutions under suitable geometric 
conditions. 
\end{abstract}

\maketitle
\thispagestyle{empty} 
\def\IMSmarkvadjust{0 pt}
\def\IMSmarkhadjust{0 pt}
\def\IMSmarkhpadding{0 pt}
\def\IMSpubltext{Published in modified form:}
\def\SBIMSMark#1#2#3{
 \font\SBF=cmss10 at 10 true pt
 \font\SBI=cmssi10 at 10 true pt
 \setbox0=\hbox{\SBF \hbox to \IMSmarkhpadding{\relax}
                Stony Brook IMS Preprint \##1}
 \setbox2=\hbox to \wd0{\hfil \SBI #2}
 \setbox4=\hbox to \wd0{\hfil \SBI #3}
 \setbox6=\hbox to \wd0{\hss
             \vbox{\hsize=\wd0 \parskip=0pt \baselineskip=10 true pt
                   \copy0 \break%
                   \copy2 \break%
                   \copy4 \break}}
 \dimen0=\ht6   \advance\dimen0 by \vsize \advance\dimen0 by 8 true pt
                \advance\dimen0 by -\pagetotal
	        \advance\dimen0 by \IMSmarkvadjust
 \dimen2=\hsize \advance\dimen2 by .25 true in
	        \advance\dimen2 by \IMSmarkhadjust

%
%
  \openin2=publishd.tex
  \ifeof2\setbox0=\hbox to 0pt{}
  \else 
     \setbox0=\hbox to 3.1 true in{
                \vbox to \ht6{\hsize=3 true in \parskip=0pt  \noindent  
                {\SBI \IMSpubltext}\hfil\break
                \input publishd.tex 
                \vfill}}
  \fi
  \closein2
  \ht0=0pt \dp0=0pt
 \ht6=0pt \dp6=0pt
 \setbox8=\vbox to \dimen0{\vfill \hbox to \dimen2{\copy0 \hss \copy6}}
 \ht8=0pt \dp8=0pt \wd8=0pt
 \copy8
 \message{*** Stony Brook IMS Preprint #1, #2. #3 ***}
}

\def\IMSmarkvadjust{-30pt}
\SBIMSMark{2003/02}{October 2003}{}

\begin{quote}
\begin{small}
\end{small}
\end{quote}

\section{Introduction}
We define and study a new dynamical system in the 
space of K\"ahler metrics that represent a fixed cohomology class of a given
closed complex manifold of K\"ahler type. The critical points of this flow
are extremal metrics, that is to say, minimizers of the functional defined
by the $L^2$-norm of the scalar curvature. We derive the equation, describe 
some of its 
general properties, and prove that given an initial data, the equation 
has a unique classical solution on some time interval. It would be of great
interest to know if the solution exist for all time, or whether it
develops some singularities in finite time. We have no general answer to 
this yet. However, we show some evidence indicating that  
in some specific cases, the solution should exist for all 
time and converge to an extremal metric as time goes to infinity.

In order to put our equation in proper perspective, we begin by
recalling a different but related one, the Ricci flow. Let $M$ be a compact 
manifold $M$ of dimension $n$. Given a metric $g$, we denote its Ricci 
tensor by $Ricci_g$ and its
average scalar curvature by $r_g$. The Ricci flow 
$$\frac{dg}{dt}=2\left( \frac{r_g}{n}g-Ricci_g\right) \, ,$$
was introduced by R. Hamilton \cite{ha} as a mechanism to improve the 
properties of its initial data. It is a non-linear heat equation in the
metric, which hopefully becomes better as time passes by in the same way as 
the heat equation improves an initial distribution of heat in a given 
region, and makes it uniform all throughout as time goes to infinity. 
Hamilton used it to show that on a three dimensional manifold, an initial
metric of positive Ricci curvature flows according to this equation 
towards a limit that has constant positive sectional curvature.

In the case of a K\"ahler manifold, Hamilton's flow equation may be used
when seeking a K\"ahler-Einstein metric on the said manifold. Of course, 
this would a priori require that the first Chern class $c_1$ be signed, so
that it may be represented by K\"ahler-Einstein metrics, or their opposites.
Regardless of that consideration, the idea inspired Cao \cite{cao} to study 
the equation
$$\frac{d\o}{dt} =\eta -\rho_t \, ,$$ 
for $\eta$ a fixed real closed $(1,1)$-form representing the class 
$c_1(M)$. Here, $\o$ and $\rho$ are the K\"ahler and Ricci forms of the metric,
respectively. Using Yau's work on the Calabi 
conjecture, he proved that solutions exists for all 
$t\geq 0$ and that the path of metrics so defined 
converges to a K\"ahler metric with prescribed Ricci form $ \eta$ as 
$t\rightarrow \infty$. He went on and, under the assumption that $c_1(M)<0$, 
replaced $\eta$ in the equation above by $-\o_t$ and proved that the 
corresponding solution to the initial value problem exists for all time
and converges to a K\"ahler-Einstein metric as $t\rightarrow \infty $, 
rederiving the now famous version of the theorem of Yau and Aubin. 

As good as the Ricci flow may be, when it converges it can only do so
to an Einstein metric, fact that is only possible under {\it a priori} 
conditions on $M$ that might not hold in general. We are interested in 
studying polarized K\"ahler manifolds, and in finding canonical 
representatives of the polarization. Our flow is adapted to accomplishing 
that other somewhat different goal.

For if $(M,J,\Omega)$ is a polarized K\"ahler manifold, for a variety of
technical reasons to be clarified below, we consider metrics that are 
invariant under the action of a maximal compact subgroup of the automorphism 
group of $(M,J)$, and introduce the equation
$$\frac{d\o}{dt}=\Pi_t \rho -\rho_t\, ,$$
with initial condition a given metric representing $\Omega$.
Here $\Pi_t$ is a metric dependent projection operator that intertwines 
the metric trace with the $L^2$-orthogonal projection $\pi_t$ onto the 
space of real holomorphy potentials, these being those real valued functions 
whose gradients are holomorphic vector fields. The projection $\Pi_t$ is 
such that $\Pi_t \rho_t - \rho_t$ is cohomologous to zero, and so all metrics 
satisfying the equation represent $\Omega$. 
As such, the critical points of the equation, metrics for which 
$\Pi \rho=\rho$, are precisely those metrics whose scalar curvatures have 
holomorphic 
gradients, or said differently, the extremal metrics of Calabi 
\cite{ca1}. This fact constitutes the guiding principle behind our 
consideration of this new flow equation. 

In general, our flow equation is different from the K\"ahler version of the
Ricci flow, even when $\Omega =\pm c_1$. 
This last assertion is illustrated, for instance, by the blow-up of 
${\mathbb C}{\mathbb P}^2$ at one or two points, and the reason is basically 
a simple one: extremal metrics, which is what we seek when we consider the 
new flow, is a concept that imposes milder conditions than those required 
for the metric to be K\"ahler-Einstein, and when $\Omega=c_1$, the two 
concepts agree only if we know of additional restrictions on $c_1$ 
\cite{fu,ss}. These do not hold on the two mentioned examples. On the other 
hand, the flows always coincide on Riemann 
surfaces (where they also agree with the 
two-dimensional version of the Yamabe flow \cite{ha2}) because 
regardless of the metric $g$ under consideration, 
the holomorphy potential $\pi_g s_g$, $s_g$ the scalar curvature of $g$, 
is a topological constant determined by the Euler characteristic.
In these cases, it is still of some interest to point out that 
the extremal flow always arises as the gradient flow of a Riemannian 
functional, statement that in this dimension cannot be ascertained for the 
other two flows in light
of the Gauss-Bonnet theorem.

The main point of the present article will be to show that solutions 
to the extremal flow equation exists locally in time. However, even if these
were going to exist globally, we should not expect that they would 
converge as time approaches infinity in all possible cases. We already 
know of examples of K\"ahlerian manifolds that do not admit extremal metrics 
\cite{bb,le,ti}.

We do have already a partial picture that
explains why these examples exist. Those in \cite{le} fail 
to satisfy a necessary condition on the space of holomorphic vector fields,
while those in \cite{bb} and \cite{ti} are related to stability of the
manifold under deformations of the complex structure, property that
is independent of and not reflected by those of the Lie algebra of 
holomorphic vector fields. 

Pointing more to the heart of the problem that interest us here, we had 
proven \cite{cs2} that the set of K\"ahler 
classes that can be represented by extremal metrics is open in the 
K\"ahler cone. The study of the extremal flow equation above, and its 
potential convergence to a limit extremal metric, can be seen as a general
method that could decide if the the extremal cone is ---or is not--- closed 
also.

In proving local time existence of the extremal flow, we also
show that if the lifespan is finite, then the point-wise norm of the 
curvature tensor must blow-up as times approaches it. We leave for later
the analysis of global existence and convergence under suitable geometric 
conditions, in particular, the analysis of these issues for surfaces with 
positive first Chern class, where our flow equation seems to be a 
promising tool for the resolution 
of the extremal metric problem.

We organize the paper as follows: in \S\ref{ekm} we recall the notion and
basic facts about extremal metrics; in \S\ref{sheat} we explain in detail
the derivation of the extremal flow equation, and prove general results 
about it; in  \S\ref{sl} we linearize this flow equation, showing that
it results into a pseudo-differential perturbation of the standard time
dependent heat equation. This form of the linearization is an essential 
fact in our
proof of local time existence, done in \S\ref{loc} via a fixed point type
of argument. We end with some remarks justifying our hope that solutions to
the extremal flow will converge to an extremal metric under suitable general 
geometric conditions.

\section{Extremal K\"ahler Metrics}\label{ekm}
Let $(M,J,g)$ be a closed K\"ahler manifold of complex dimension $n$. This 
means that $(M,J)$ is a closed complex manifold and that 
$\omega (X,Y) := g(JX,Y)$, which is skew-symmetric because $g$ is a 
{\it Hermitian} Riemannian metric, is a closed 2-form. The form $\omega$ is 
called the {\it K\"ahler form}, and its  cohomology class 
$[\omega ] \in H^{2}(M,{\mathbb R})$ is called the {\it  K\"{a}hler class}. 

We denote by ${\mathfrak h}={\mathfrak h}(M,J)$ the complex Lie algebra 
of holomorphic vector fields of $(M,J)$. Since $M$ is compact, this is
precisely the Lie algebra of the group of biholomorphism of $(M,J)$.
The subset ${\mathfrak h}_0$ of holomorphic vector fields with zeroes is an 
ideal of ${\mathfrak h}$, and the quotient algebra 
${\mathfrak h}/{\mathfrak h}_0$ is Abelian. 

We shall 
say that $(M,J)$ is a {\it generic} complex manifold of K\"ahler type 
if ${\mathfrak h}_0$ is trivial. Typically, complex manifolds 
carry no non-zero holomorphic vector fields and are, therefore, generic in our 
sense. However, our definition includes also those complex manifolds whose
non-zero holomorphic fields have empty zero sets.

By complex  multi-linearity, we may extend the metric
 $g$, the Levi-Civita connection $\nabla$ and the curvature 
tensor $R$  to the complexified tangent bundle
${\mathbb C} \otimes TM$. Since 
 ${\mathbb C} \otimes TM$ decomposes
into the $\pm i$-eigenspaces of $J$,  ${\mathbb C}
 \otimes TM = T^{1,0}M \oplus T^{0,1}M$, 
we can express any tensor field or differential operator in terms of 
the corresponding decomposition.
For example, if $\{z^{1}, \ldots ,z^{n}\}$ is a holomorphic coordinate system 
on $M$, we get induced bases  $\{ \d_{z^\j} \}$ and $
\{ \d_{z^{\b{\j}}} :=\d_{\b{z}^\j} \}$ for $T^{1,0}M$ 
and $T^{0,1}M$, respectively, and if we express the metric  $g$ in terms of
this basis by setting $g_{\mu\nu} :=g\left( \d_{z^{\mu}},
\partial _{z ^{\nu}}\right)$, 
where the indices 
$\mu$, $\nu$ range over  $\{ 1, \ldots ,n,\bar{1}, \ldots ,\bar{n}\}$,
it follows from  the Hermiticity condition 
 that $g_{\j k} = g_{\b{\j}\b{k}}=0$, and that
$\omega = \omega _{{j}\b{k}}dz^{j}\wedge d\b{z}^{k}= i
g_{\j \b{k}}dz^{\j}\wedge d\b{z}^{k}$. 

The complexification of the exterior algebra can be decomposed into
a direct sum of of forms of type $(p,q)$. Indeed, we have
$\wedge ^r M =\bigoplus_{p+q=r} \wedge ^{p,q}M$. The integrability of
$J$ implies that the exterior derivative  
 $d$ splits as $d={\partial}+\bar{\partial}$,
where ${\partial}: \wedge ^{p,q}M\to \wedge ^{p+1,q}M$,
$\bar{\partial}: \wedge ^{p,q}M \to \wedge ^{p,q+1}M$,
${\partial}^2=\bar{\partial}^2=0$ and ${\partial}\bar{\partial}=
-\bar{\partial}{\partial}$. Complex conjugation also extends, and we define
a form to be real if it is invariant under this operation.
An important result in K\"ahler geometry is that,
given a $d$-exact real form $\beta $ of type $(p,p)$, there exists a real form
$\a $ of type $(p-1,p-1)$ such that $\beta = i\ddb \a $.

The {\it Ricci form} $\rho$ is defined in terms of the 
Ricci tensor $r$ of  $g$ by $\rho (X,Y) =  r (JX,Y)$. It is a closed form 
whose components are given by  
$$r_{j\overline{k}} = -i \rho_{j\overline{k}} = 
-\frac{\partial ^{2}}{\partial {z^{j}} \partial 
{\overline{z}^{k}}} \log{{\rm  det}(g_{p\overline{q}})}\, .$$ 
The form ${\displaystyle \frac{1}{2\pi}}\rho$ is the curvature of the 
canonical line bundle $\kappa= \Lambda^n (T^*M)^{1,0}$,
and represents the first Chern class $c_1=c_1(M,J)$.

The {\it scalar curvature} $s$ is,  by definition,  the trace
$s = r^{\mu}_{\mu}= 2g^{\j\b{k}}r_{\j\b{k}}$ 
of the Ricci tensor, and can be conveniently calculated by the
formula 
\begin{equation}
s \, \omega^{\w n} = 2n \, \rho \w \omega^{\w (n-1)}\, .\label{eq:sca}
\end{equation}
Since the volume form is given by $d\mu ={\displaystyle \frac{\omega^{\w n}}
{n!}}$ and $M$ is compact, this formula implies that
\begin{equation}
 \int_M s \, d\mu = \frac{4\pi}{(n-1)!}  c_1
 \cup [\omega]^{\cup (n-1)}\, , \label{gb}
\end{equation}
a quantity that only depends upon the complex structure $J$ and the 
cohomology class $[\o]$, and that generalizes the well-known Gauss-Bonnet
theorem for surfaces. Notice that ${\displaystyle \int_M d\mu =
 \frac{1}{n!} [\omega]^{\cup n}}$, and so {\it the average scalar curvature}
\begin{equation}
s_0= 4\pi n \frac{c_1 \cup [\omega]^{n-1}}{[\o]^{n}}\, ,\label{prs}
\end{equation}
is also a quantity that depends only on the K\"ahler class $[\omega]$ 
and the homotopy class of the complex-structure tensor $J$.

Suppose that  $(M,J)$ is {\it polarized}
by a positive class $\Omega \in H^{1,1}(M, {\mathbb C})\cap  
H^{2}(M, {\mathbb R})$. Let ${\mathfrak M}_{\Omega}$ be the set of all 
K\"ahler forms representing $\Omega$. Since any two elements $\tilde{\o}$ and
$\o$ of ${\mathfrak M}_{\Omega}$ are such that $\tilde{\o}=\o+i\ddb \varphi$
for some real valued potential function $\v$, at the expense of 
fixing a background K\"ahler metric $\omega$ that represents $\Omega$, we 
can describe ${\mathfrak M}_{\Omega}$ as ${\mathfrak M}_{\Omega}=\{ 
\o_{\varphi}=\o +i\ddb \varphi: \; \o_{\varphi}>0
\}$. Thus, ${\mathfrak M}_{\Omega}$ is an affine space modeled on
an open subset of the space of smooth functions that parametrizes the 
deformations of the base point $\o$. We may topologize it 
by defining a suitable topology on the space of deformation potentials.

In what follows, we shall not distinguish between the 
K\"ahler metric and its K\"ahler form, passing from one to the other at
will.

Consider the functional  
\begin{equation}
 \begin{array}{rcl}
{\mathfrak M}_{\Omega} & \stackrel{\Phi_{\Omega}}{\longrightarrow} & 
{\mathbb R} \vspace{.08in} \\
\o & \mapsto & {\displaystyle \int _M s_\o ^2 d{\mu}_\o } \end{array}\, ,
\label{ene}
\end{equation}
where the metric associated with the form $\o$
has scalar curvature $s_\o$  and  volume form $d{\mu}_\o$. A critical
point of this functional is by definition
an {\it extremal} K\"ahler metric \cite{ca1}. They were introduced by E. Calabi
with the idea of seeking canonical representatives of $\Omega$.

Given any K\"ahler metric $g$, a smooth complex-valued
function $f$ gives rise to the (1,0) vector field  $f\mapsto \partial^{\#} 
f=\partial^{\#}_g f$ defined by the expression
$$g(\d^{\#} f, \hspace{1mm}\cdot \hspace{1mm})=\db f \, .$$
 This vector field is {\it holomorphic} iff we require 
that $\b{\partial}\partial^{\#} f=0$, condition equivalent 
to $f$ being in the kernel of the operator
\begin{equation}
L_g f:= (\b{\d}{\d}^{\#})^{\ast}\b{\d}{\d}^{\#}f= \frac{1}{4}\lap^2 f + 
\frac{1}{2}r ^{\mu\nu}\n_{\mu}\n_{\nu}f+ \frac{1}{2}(\n^{\b{\ell}}\s )
\n_{\b{\ell}}f \, .\label{lic}
\end{equation}
The ideal ${\mathfrak h}_0$ consists of vector fields of the
form $\partial^{\#}_g f$, for a function $f$ in the kernel of $L_g$. 
In other words, a holomorphic vector field $\Xi$ can be written
as $\partial_g^{\#}f$ iff the zero set of $\Xi$ is non-empty.

The variation of $\Phi_{\O}$ can be given in terms of the operator
$L_g$ above. Indeed, 
$$\frac{d}{dt}\Phi_{\Omega}(\o+ti\ddb \v)\mid_{t=0}=-4\int s_{\o}L_{\o}\v 
d\mu_{\o} \, .$$
Hence, the Euler-Lagrange's equation for a critical point $g$ of (\ref{ene}) 
can be cast as the fact that the scalar curvature $s_g$ is a real valued 
function in the kernel of $L_g$. Thus, the vector field $\d^{\#}_g s_g$ must 
be holomorphic and its imaginary part is a Killing field of $g$.

Given a K\"ahler metric $\o$, its normalized Ricci potential $\psi_\o$ is 
defined to be the only function orthogonal to the constants such that
$\rho=\rho_{H}+i\ddb \psi_{\o}$, where $\rho_{H}$ is the 
$\o$-harmonic component of $\rho$. In terms of the scalar curvature and 
its projection onto the constants, we have that 
$\psi_{\o}=-G_{\o}(s_{\o} -s_0)$. If ${\mathcal K}$ denotes the 
K\"ahler cone of $(M,J)$, the {\it Futaki character} \cite{fu} 
is defined to be the map 
\begin{equation}
\begin{array}{c}
{\mathfrak F}: {\mathfrak h} \times {\mathcal K} \longrightarrow   
{\mathbb C} \\ {\mathfrak F}(\Xi,[\o])=
2{\displaystyle \int _{M}
\Xi(\psi _{\o})d\mu = -2\int _{M} \Xi( G_\o (s_\o-s_0) ) d\mu_\o }\, .
\end{array}
\label{futc}
\end{equation}
It is calculated using a particular representative $\o$ of the class, but
it depends only on $\O$ and not on the particular choice of representative
\cite{fu,ca2}. And when applied to a 
holomorphic vector field of the form $\Xi = \d ^{\#}f$, it yields 
\begin{equation}
{\mathfrak F} (\Xi,[\o ]) = {\displaystyle -\int_Mf(s_\o-s_0) \, d\mu_\o }
\, .\label{fuss}
\end{equation}

If a metric $\o\in {\mathfrak M}_{\O}$ is extremal, we may apply (\ref{fuss}) 
to the vector field $\d^{\#}_{\o}s_{\o}$ and obtain that
${\mathfrak F} (\d^{\#}_{\o}s_{\o},\O)=-\| s_{\o}-s_0\|^2$. Thus, the 
Futaki character represents an obstruction to $\o$ being a metric of
constant scalar curvature.

Extremal metrics achieve the infimum of $\Phi_{\Omega}$ over 
${\mathfrak M}_{\Omega}$. In fact, the critical value $E(\O)$ they achieve 
is a differentiable function of $\O$ \cite{si3}, and there exists 
a holomorphic vector field $X_{\Omega}$ \cite{fm} such that 
\begin{equation}
\Phi_{\Omega}(\o)\geq E(\Omega):=
s _{0}^2 \frac {\Omega^n}{n!}-
{\mathfrak F} (X_{\Omega},\Omega) 
\label{lb}
\end{equation}
for all $\o \in {\mathfrak M}_\Omega$. 

\section{Derivation of the evolution equation} \label{sheat}
Calabi \cite{ca2} showed that the identity component of the isometry group of
an extremal K\"ahler metric $g$ is a maximal compact subgroup of the
identity component of the biholomorphism group of  $(M,J)$. This 
implies that, up to conjugation, the identity components of the isometry
groups of extremal K\"ahler metrics coincide \cite{cs2}. Therefore, 
modulo biholomorphisms, the search for extremal K\"ahler metrics is 
completely equivalent to the search for extremal metrics among those
that are invariant under the action of a fixed maximal compact subgroup
of the connected biholomorphism group. This last problem, however, turns out
to be technically easier to analyze.

\subsection{Real holomorphy potentials}
Given any K\"ahler metric $g$ on $(M,J)$, every complex-valued function 
$f$ in the kernel of $L_g=(\b{\d}{\d}^{\#}_g)^{\ast}\b{\d}{\d}^{\#}_g$ in 
(\ref{lic}) is associated with the holomorphic vector field 
$\Xi= {\partial}^{\#}_gf$, and since the operator $L_g$ is elliptic,
the space of such functions is finite dimensional. However, 
$L_g$ is not, generally  speaking, a real operator. Therefore, the real
and imaginary parts of a function in its kernel do not have to be elements
of the kernel also. It has been proven elsewhere \cite{cs2} that if $f$ is a 
real valued function in the kernel of $L_g$, 
then the imaginary part of $\partial^{\#}_g f$ is a Killing field of
$g$, and that a Killing field arises in this way if, and only if, its
zero set is not empty.

Let $G$ be a maximal compact subgroup of the biholomorphism group of
$(M,J)$, and $g$ be a K\"ahler metric on $M$ with K\"ahler class $\Omega$.
Without loss of generality, we assume that $g$ is $G$-invariant. 
We denote by $L^2_{k,G}$ the real Hilbert space of $G$-invariant 
real-valued functions of class $L^2_{k}$, and consider $G$-invariant 
deformations of this metric preserving the K\"ahler class:
 \begin{equation}
\tilde{\omega }=\o +i\ddb \varphi \, , \quad \varphi \in L_{k+4,G}^2\, , \;
k>n.
\label{met}\end{equation}
In this expression, the condition $k>n$ ensures that $L_{k,G}^2$ is a Banach
algebra, making the scalar curvature of 
$\tilde{\o}$ a well-defined function in the space. 

We denote by $\mathfrak g\subset {\mathfrak h}$ the Lie algebra of $G$, and by 
${\mathfrak z}$ the center of ${\mathfrak g}$.
We let ${\mathfrak z}_0= {\mathfrak z}\cap {\mathfrak g}_0$,
where ${\mathfrak g}_0\subset {\mathfrak g}$ is the ideal of 
Killing  fields which have  zeroes. 
If $\tilde{g}$ is any $G$-invariant K\"ahler metric on $(M,J)$, then 
each element
of ${\mathfrak z}_0$ is of the form $J\,\nabla_{\tilde{g}}f$ for a real-valued
solution of (\ref{lic}). In fact, 
${\mathfrak z}_0$ corresponds  to the set of real solutions
$f$ which are {\em invariant under $G$}, since   
$${\partial}^{\#}: \ker
 [(\b{\d}{\d}_{\tilde{g}}^{\#})^{\ast}\b{\d}{\d}_{\tilde{g}}] \to 
{\mathfrak h}_0$$
is a homomorphism of $G$-modules.

The restriction of $\ker (\b{\d}{\d}_{\tilde{g}}^{\#})^{\ast}\b{\d}
{\d}_{\tilde{g}}$
to $L^2_{k+4,G}$ depends smoothly on the $G$-invariant metric 
$\tilde{g}$. Indeed, 
choose  a basis $\{ X_1, \ldots, X_m\}$ for ${\mathfrak z}_0$,
and, for each  $(1,1)$-form $\chi$ on $(M,J)$, consider the set 
of functions 
$$\begin{array}{rcl}
p_0(\chi ) & = &  1 \\ p_j ( \chi ) & = &
2i G_g{\db}^{*}_{g}((JX_{j}+iX_{j})\hok \chi )
\, , \quad j=1, \ldots , m
\end{array}$$ where $G_g$ is the Green's operator of the
metric $g$. 
If $\tilde{\o }$  is the K\"ahler form of the $G$-invariant metric 
$\tilde{g}$, then $\d ^{\#}_{\tilde{g}}p_j(\tilde{\o})=JX_{j}
 +iX_{\jmath}$,
and the set $\{ p_j(\tilde{\o })\} _{j=0}^{m}$ consists of real-valued 
functions and forms a basis for 
$\ker (\b{\d}{\d}_{\tilde{g}}^{\#})^{\ast}\b{\d}
{\d}_{\tilde{g}}$. Furthermore, for metrics $\tilde{\o}$ as in (\ref{met}),
the map $\varphi \mapsto p_j(\o +i \ddb \varphi)$
is, for each $j$, bounded as a linear map from $L^2_{k+4,G}$ to 
$L^2_{k+3, G}$.

With respect to the  fixed $L^2$ inner product, 
let $\{ f_{\tilde{\o }}^0, \ldots , f_{\tilde{\o }}^m\}$ be the 
orthonormal set extracted from $\{ p_j (\tilde{\o }) \}$ by 
the Gram-Schmidt procedure. We then let
\begin{eqnarray} \pi_{\tilde{\o }}: L^2_{k,G}&\to &L^2_{k,G}\nonumber \\
u&\mapsto & \sum_{j=0}^m \<  f_{\tilde{\o }}^j, u\>_{L^2} f_{\tilde{\o }}^j
\label{proj1}\end{eqnarray}
denote the associated projector. In fact, by the regularity of the 
functions $\{ p_1, \ldots , p_{m}\}$, this projection can be defined on
$L^2_{k+j,G}$ for $j=0,\, 1,\, 2,\, 3$, and for metrics as in 
(\ref{met}), the map $\varphi \mapsto \pi_{\tilde{\o }}$ is smooth 
from a suitable neighborhood
of the origin in $L^2_{k+4,G}$ to the real
Hilbert space $\mbox{End} (L^2_{k+j,G})\cong \bigotimes^2 L^2_{k+j,G}$. 

The holomorphic vector field $X_{\O}$ of the class that yields the lower bound
(\ref{lb}) is given by $X_{\Omega}=\d_g^{\#}(\pi_g s_g-s_0)$.
As such, it may depend on the choice of maximal compact subgroup $G$ of the 
automorphism group of $(M,J)$, but the value of 
${\mathfrak F} (X_{\Omega},\Omega)$ does not. The critical value
$E(\O)$, or energy of the class $\O$, is nothing more than 
$$E({\Omega})= \int (\pi_g s_g)^2 d\mu_g \, .$$
(This way of computing the energy of the class through $G$-invariant metrics is
very convenient and has been used several times elsewhere \cite{si2,ss,ss2} 
with other (but related) purposes in mind.)

\begin{remark}
For generic manifolds $(M,J)$, the space of real holomorphy potentials consists
of the constant functions only. Thus, $\pi_g s_g$ is the
constant given in (\ref{prs}) no matter what $G$-invariant representative
$g$ of the class $\O$ we consider. However, this last function could be
constant even under the presence of non-trivial holomorphic vector fields 
with zeroes. For that, the Futaki invariant of the class must vanish, and in
that case, the vector field $X_{\Omega}$ associated to $\O$ is trivial.

A simple instance of this is given by a compact Riemann surface 
$\Sigma$, where for any $G$-invariant metric $g$, we have that
$\pi_g s_g =4\pi \chi(\Sigma)/\mu_{g}(\Sigma)$ where $\chi(\Sigma)$ is the
Euler characteristic. Indeed, Riemann surfaces are either hyperbolic, parabolic
or elliptic. The first two of these are generic, and the assertion is clear 
then. For the remaining case, that of the Riemann sphere, the space 
${\mathfrak h}_0$ is non-trivial but the Futaki character vanishes. The
assertion follows from the Gauss-Bonnet theorem embodied in (\ref{gb}), and 
the identity (\ref{fuss}) applied to $f=\pi_g s_g$.
\end{remark}
\medskip

From now on, we shall denote by ${\mathfrak M}_{\Omega,G}$ the set of 
$G$-invariant K\"ahler metrics representing the class $\Omega$. Given a 
path of metrics $\o_t \in {\mathfrak M}_{\Omega, G}$, since the 
kernel of $(\b{\d}{\d}_{\tilde{g}}^{\#})^{\ast}\b{\d}
{\d}_{\tilde{g}}$ depends smoothly on $\tilde{g}$, the differential of 
$\pi_t s_t$ is well-defined. Here, $\pi_t$ and $s_t$ are the projection 
operator and scalar curvature associated to the metric $\o_t$, 
respectively. Since $\pi_t s_t$ is of order four in the potential of the 
metric, naively we would expect its differential to be an operator of order 
four on the tangent space to ${\mathfrak M}_{\Omega, G}$ at $\o_t$. However, 
we get something significantly better, and gain quite a bit of regularity.
This fact will be very convenient later on.

\begin{lemma}
\label{le1}
Let $\o_t=\o+i\ddb \v_t $ be a path of metrics in 
${\mathfrak M}_{\Omega, G}$ with $\o_0=\o$. Consider the projection 
$\pi_t s_t$ of the scalar curvature $s_t$ onto the space of real 
holomorphy potentials, and let
$\dot{\v}_t = \frac{d}{dt}\v_t$. Then
$$\frac{d}{dt}(\pi_t s_t) = \d \dot{\v}_t \hok X_{\Omega}=(\d^{\#}
\dot{\v}_t, X_{\Omega})_t=(\d \dot{\v}_t, \d (\pi_t s_t) )_t\, ,$$
where $X_{\Omega}=\d_t^{\#} \pi_t s_t$ is the holomorphic vector field
of the class $\Omega$.
In particular, this derivative is a differential operator of 
order one in $\dot{\v}_t$ whose coefficients depend non-linearly on the 
metric $\o_t$. 
\end{lemma}

{\it Proof}. By the invariance of the Futaki character, if 
$\pi_\o s_\o$ is constant then so will be $\pi_{\tilde{\o}}s_{\tilde{\o}}$
for any other metric $\tilde{\o}$ in ${\mathfrak M}_{\Omega, G}$ (see
\S4 of \cite{ss}). In that case, $X_{\Omega}$ is trivial
and both sides of the expression in the statement are zero. The result 
follows.

So let us assume that $\pi_t s_t$ is not constant. Here and below 
we use the subscript $t$ to denote geometric
quantities associated with $\o_t$. Thus, the imaginary part
of $X_{\Omega}=\d_{\o_t}^{\#} \pi_t s_t $ is a non-trivial Killing vector 
field, and in the construction of the projection map above, we can choose
a basis $\{ X_{j}\}$ for ${\mathfrak z}_0$ such that 
$X_{\Omega}=\d_{\o_t}^{\#}(\pi_t s_t) = JX_1 +iX_1=X_{\Omega}$. Hence,
$$\pi_t s_t = 2i G_{t}\db_{t}^{*}(\o_t \hok X_{\Omega})+s_0\, ,$$
where $s_0$ is the projection (\ref{prs}) of $s$ onto the constants, a 
function that depends on $\O$ and $J$ but not on the particular choice of 
metric in ${\mathfrak M}_{\O,G}$.

By the K\"ahler identity $\db_{t}^{*}=-i[\Lambda_{t},\d]$, we conclude that
$$\pi_t s_t = 2G_{t}\Lambda_t \d (\o_t \hok X_{\Omega})+s_0$$
and, therefore,
$$\frac{d}{dt}\pi_t s_t = 2G_{t}\Lambda_t \d (\dot{\o}_t \hok X_{\Omega})
+2G_t \dot{\Lambda}_t \d  (\o_t \hok X_{\Omega}) +
2\dot{G}_t \Lambda_t \d (\o_t \hok X_{\Omega})\, .$$

The last two terms in the expression above cancel each other out. 
Indeed, $\o_t \hok X_{\Omega}=-i\db(\pi_t s_t)$ and computing the derivative
of $\dot{\Lambda}_t$ in terms of $\dot{\v}_t$, we see that
$2G_t \dot{\Lambda}_t \d  (\o_t \hok X_{\Omega})=2G_t(i\ddb \dot{\v}_t,
i\ddb (\pi_t s_t))_t$. On the other hand, the differential of the Green's 
operator is given by $-G_t \dot{\Delta}_t G_t$, and we obtain that
$2\dot{G}_t \Lambda_t \d (\o_t \hok X_{\Omega})=
2G_t\dot{\Delta}_t G_t \Lambda_t i\ddb (\pi_t s_t) = -G_t\dot{\Delta}_t(\pi_t 
s_t) =-2G_t(i\ddb \dot{\v}_t,i\ddb (\pi_t s_t))_t$.

Since the real and imaginary parts of $X_{\Omega}$ are Killing vector fields
and the metric potential $\v_t$ is $G$-invariant, we have that 
$X_{\Omega}(\dot{\v}_t)=0$, and so $\d \dot{\v}_t \hok 
X_{\Omega}=(\d \dot{\v}_t,i\d (\pi_t s_t))_t$ is orthogonal to the 
constants. On the other hand, since $X_{\Omega}$ is holomorphic, we have that
$\dot{\o}_t\hok X_{\Omega}=i\ddb \dot{\v}_t \hok X_{\Omega}= -i 
\db(\d \dot{\v}_t \hok X_{\Omega})$. Hence, 
 $2G_{t}\Lambda_t \d (\dot{\o}_t \hok X_{\Omega})=-2G_t \Lambda_t i\ddb 
(\d \dot{\v}_t \hok X_{\Omega})$, and the desired result follows now because
$G_t$ is the inverse of the Laplacian in the complement of the constants. \qed

Given any K\"ahler metric $g$ in ${\mathfrak M}_{\Omega,G}$, we have seen that
the extremal vector field $X_{\Omega}$ of the class can be written as 
$X_{\Omega}=\d_g^{\#} (\pi_g s_g)$. Thus, the critical points
of $\pi_g s_g$ corresponds to zeroes of $X_{\Omega}$, and are therefore,
independent of $g$. We may use the Lemma above to strengthen this assertion
a bit, and derive the following remarkable consequence. This result is
reminiscent of the convexity theorem on the image of moment mappings
\cite{at,gust}.

\begin{theorem} \label{tc}
Let $\o$ be any metric in ${\mathfrak M}_{\Omega, G}$ and consider the
function $\pi_{\o} s_{\o}$ obtained by projection of the scalar 
curvature onto the space of real holomorphy potentials. Then the range
of $\pi_\o s_\o$ is a closed interval on the real line that only
depends on the class $\Omega$ and the complex structure $J$, but not on the 
particular choice of metric $\o \in {\mathfrak M}_{\Omega, G}$. 
\end{theorem}

{\it Proof}. Let $\o_t=\o+i\ddb \v_t$ be a path in 
${\mathfrak M}_{\Omega, G}$. By Lemma \ref{le1}, we have that
$$\frac{d}{dt}\pi_t s_t = (\d \dot{\v}_t , \d (\pi_t s_t) )_t \, .$$
Since the maximum and minimum of $\pi_t s_t $ occur at critical points, this
expression shows that these extrema values do not change with $t$. The
result follows because ${\mathfrak M}_{\Omega, G}$ is path connected. \qed

The projection $\pi_g$ onto holomorphy potentials has an appropriate lift
$\Pi_g$ to the level of $G$-invariant $(1,1)$ forms, which we discuss now.
We denote by $\wedge ^{1,1}_{k,G}$ the space of real forms 
of type $(1,1)$ that are invariant under $G$ and of class $L_{k}^2$. 

\begin{lemma}
\label{le2}
Given any $G$-invariant metric $\tilde{\o}$, there exists a uniquely defined
continuous projection map 
\begin{equation}
\Pi _{\tilde{\o}}:\wedge ^{1,1}_{k+2,G} \mapsto 
\wedge ^{1,1}_{k+2,G} \, , \label{proj} \end{equation}
that intertwines the trace with the projection map $\pi _{\tilde{\o}}$ in
{\rm (\ref{proj1})}, and it is  
such that $\eta -\Pi _{\tilde{\o}}\eta $ is cohomologous to zero for all 
$\eta \in  \wedge ^{1,1}_{k+2,G}$. For metrics $\tilde{\o}$ as in
{\rm (\ref{met})}, the map $\varphi \mapsto \Pi _{\tilde{\o}}$ from
$L^2_{k+4,G}$ to ${\rm End} (\wedge ^{1,1}_{k+2,G})$ is smooth. 
\end{lemma}

{\it Proof}. Let $\eta \in \wedge ^{1,1}_{k+2,G}$. Since 
$\Pi _{\tilde
{\o}}\, \eta $ must be of the form $\eta + i\ddb f$ for some real valued 
function $f$, the intertwining property of the projection and trace gives that
$${\rm trace}_{\tilde{\o}}\, \eta -\frac{1}{2}\Delta _{\tilde{\o}}f=
\pi _{\tilde{\o}}\,  {\rm trace}_{\tilde{\o}}\, \eta \, ,$$ and so
$$\Delta _{\tilde{\o}}f=-2(\pi_{\tilde{\o}} -1){\rm trace}_{\tilde{\o}}\eta 
\, .$$ 
The right side of this expression is a $G$-invariant real valued function
in the complement of the constants. Thus, we can solve the equation for $f$
to obtain a real valued function that is invariant under $G$. By the
continuity properties of the map $\pi _{\tilde{\omega}}$ for metrics as in
(\ref{met}), we conclude that $\varphi \mapsto \Pi _{\tilde{\o}}$  
is a smooth map from a suitable neighborhood of the origin in 
$L^2_{k+4,G}$ to the real Hilbert space ${\rm End} (
\wedge ^{1,1}_{k+2,G})$. \qed

\subsection{Extremal flow equation}
The projection operators $\pi_g$ and $\Pi_g$ are essential elements in our 
study of extremal K\"ahler metrics. To begin with, they lead to alternative
characterizations of the extremality condition of a metric that are quite
suitable for analytical purposes. Indeed, if $g$ is an extremal metric, it
must be invariant under a maximal compact subgroup $G$ of ${\rm Aut}(M,J)$ and
we have that
\begin{equation}
\rho_g= \Pi_g \rho_g\, .\label{re}
\end{equation}
Conversely, any metric $g$ that is $G$-invariant and satisfies this equation
must be extremal. Evidently, this tensorial equation for the $G$-invariant 
metric $g$ can be recast in terms of the equivalent scalar equation
\begin{equation}
s_g= \pi_g s_g\, ,\label{se}
\end{equation}
which also serves as a characterization of the extremal condition of $g$.
 
Notwithstanding the richer set-up, the definition of extremality embodied 
by (\ref{re}) is analogous to the definition of an Einstein metric, and 
its algebraic nature is quite suitable and direct for analytical purposes.
For instance, rephrasing the proof for Einstein metrics, we may use it
to show that any $G$-invariant metric that is a $C^{1,\alpha}$ weak 
solution of the extremal equation in harmonic coordinates, must be in
fact a smooth extremal K\"ahler metric (compare this statement with Lemma 1 
in \cite{cs0}).

However, the key use we shall make of these projections will be to define the
{\it extremal K\"ahler flow}, a tool intended to govern the improvement of
an initial representative of the class $\O$ towards one that is extremal.

The idea of using {\it good} flows to better geometric quantities was 
originally used by Eells and Sampson \cite{es} in another context, and
reconsidered by Hamilton \cite{ha} in his definition of the Ricci flow.
In our case, we are given a metric in ${\mathfrak M}_{\Omega,G}$ and try to 
improve it by 
means of a non-linear {\it pseudo-differential} heat equation, requiring 
the velocity of the curve to equal the component of the Ricci
curvature that is perpendicular to the image of $\Pi$. 

More precisely, we fix a maximal compact subgroup $G$ of the automorphism 
group of $(M,J)$, and work on ${\mathfrak M}_{\Omega,G}$, the space of all 
$G$-invariant K\"ahler forms that represent $\Omega$. Given $\o \in 
{\mathfrak M}_{\Omega,G}$, we consider a path $\o_t$ of K\"ahler 
metrics that starts at $\o$ when $t=0$ and obeys the flow equation
$\d _t \o _t  =  -\rho_t+\Pi _t\rho_t$. 
Since $-\rho_t+\Pi _t\rho_t$ is cohomologous to zero and $G$-invariant, 
for as long as the solution exists, we will have that $\o_t \in 
{\mathfrak M}_{\Omega,G}$. Thus, our evolution equation is given by
the initial value problem 
\begin{equation}
\label{evol}
\begin{array}{rcl}
\d _t \o _t & = & -\rho_t+\Pi _t\rho_t \, , \\
\o _{0} & = & \o \, .\end{array}
\end{equation}
Critical points of this equation correspond precisely to extremal metrics, 
those that satisfy (\ref{re}).

In the same manner as (\ref{re}) has the alternative scalar description
(\ref{se}), we may reformulate (\ref{evol}) as an scalar equation. If
$\o_t = \o +i\ddb \v_t$, we have that
$\Pi _t\rho_t -\rho_t=i \ddb G_t( s_t - \pi_t s_t )$, where
$G_t$ is the Green's operator of the metric $\o_t$. By compactness
of $M$, we see that the deformation potential $\v_t$ evolves according to
\begin{equation}
\label{evol2}
\begin{array}{rcl}
\d _t \varphi _t & = & G_t( s_t -\pi _t s_t ) \, , \\
\varphi_{0} & = & 0 \, . \end{array}
\end{equation}

A critical point of (\ref{evol2}) is given by a metric for which 
$G_\o (s_\o-\pi_\o s_\o)=0$. Since $s_\o-\pi_\o s_\o$ is orthogonal to the 
constant, this condition is precisely the extremality condition
(\ref{se}).

\subsection{General properties of the extremal flow}
We begin by making a rather expected observation.

\begin{proposition}
Let $\o_t$ be a solution of the initial value problem 
{\rm (\ref{evol})}. If $d\mu_t$ is the volume form, we have that
$$\frac{d}{dt} d\mu_t = \frac{1}{2}(\pi _t s_t -s_t)d\mu_t \, .$$
In particular, the volume of $\o_t$ is constant.
\end{proposition}
 
{\it Proof}. The volume form is given by
$$d\mu_t =\frac{\o_t^n}{n!}\, .$$
Differentiating with respect to $t$, we obtain:
$$\frac{d}{dt}d\mu_t = \frac{1}{(n-1)!}\o_t^{n-1}\wedge \dot{\o_t}=
\frac{1}{(n-1)!}\o_t^{n-1}\wedge (\Pi_t \rho_t-\rho_t)=
\frac{1}{2}(\pi _t s_t -s_t)d\mu_t \, ,$$
as desired. Notice that this form of maximal rank is exact. \qed

Our next results address the plausible existence of fixed points or 
periodic solutions of the flow equation. 

Observe that (\ref{evol}) is invariant under the group of diffeomorphism 
that preserve the complex structure $J$. An {\it extremal soliton} is a 
solution that changes only by such a diffeomorphism. Then, 
there must be a holomorphic vector field $V=(V^i)$ such that $V_{i,\bar{j}}+
V_{\bar{j},i}=\Pi \rho _{i\bar{j}}-\rho_{i\bar{j}}$. If the vector field
$V$ has a holomorphy potential $f$, we refer to the pair $(g,V)$ as 
a gradient extremal soliton.

\begin{proposition}
There are no extremal solitons other than extremal metrics.
\end{proposition}

{\it Proof}. Suppose we have an extremal gradient soliton $(g,V)$ 
defined by a 
holomorphic potential $f$. Then 
$$i\ddb f = \Pi \rho - \rho \, , $$
and therefore, 
$$f = G_g(s - \pi s) \, .$$
This implies that $\Delta f=s-\pi s$ and since $\Delta $ is a real operator, 
the holomorphy potential $f$ must be real. But $f$ is a holomorphy potential,
so it is $L^2$-orthogonal to $s-\pi s$. Hence,
$$\| \nabla f \| ^2 = \int f \Delta f d\mu_g = \int f(1-\pi)s\, d\mu=0 \, .$$
Thus, $f$ is constant, and therefore, necessarily zero.

Thus, a non-trivial soliton, if any, must be given by a holomorphic vector 
field $V$ that is not a gradient. The set of all such vector fields forms an 
Abelian 
subalgebra of the algebra of holomorphic vector fields. The group of 
diffeomorphism they generate must be in the maximal compact subgroup $G$ of
isometries of the metric. This vector field does not change the metric
and so $\o_t =({\rm exp}(tV))^{*}\o=\o$. Hence, $\dot{\o}_t=0=
\Pi \rho - \rho $, and the metric is extremal. \qed

\begin{remark} 
Evidently, the compactness of $M$ plays an important r\^ole in this argument
that rules out extremal solitons other than extremal metrics. They do exist 
in the non-compact case, where they give rise to certain points in the moduli 
space of these metrics.
\end{remark}

We now show that the evolution equation (\ref{evol}) is {\it essentially} 
the gradient flow of the $K$-energy, function that also serve to characterize 
extremal K\"ahler metrics \cite{si1}. Indeed, given two elements $\omega_0$ 
and $\omega_1$ of ${\mathfrak M}_{\Omega,G}$,
there exists a $G$-invariant function $\varphi$, unique modulo constants, 
such that $\omega_1 = \omega_0 +i{\partial  \overline{\partial}}\varphi$.
Let $\v_t$ be a curve of $G$-invariant functions such that 
$\omega_t=\omega_0 +i{\partial  \overline{\partial}}\v_t \in 
{\mathfrak M}_{\Omega,G}$ and $\omega(0)=\omega_0$, $\omega(1)=\omega_1$. 
We set
$$M(\omega_{0},\omega_{1})= -\int _{0}^{1}\!
dt \int _{M} \! \dot{\v}_t(s_t-\pi_t s_t)d\mu_t\, ,$$
where $s_t$ and $d\mu_t$ are the scalar curvature and volume form of 
the metric $\omega_t$, $\pi_t$ is the projection (\ref{proj1}) onto the 
space of $G$-invariant holomorphic potentials associated with this same 
metric, 
and $\dot{\varphi}_t={\displaystyle \frac{d\varphi_t}{dt}}$. This definition
is independent of the curve $t\rightarrow \v_t$ chosen. 

Fix $\o_0\in {\mathfrak M}_{\Omega,G}$. The $K$-energy 
is defined to be
\begin{equation}
\begin{array}{rcl}
{\mathfrak M}_{\Omega,G} & \stackrel{\kappa}{\longrightarrow} & 
{\mathbb R} \\
\omega & \rightarrow & M(\omega_0, \omega)\, .
\end{array}
\end{equation}
We have (see Proposition 2 in \cite{si1}) that
$$\frac{d}{dt}\kappa (\omega_t)=-\int_{M} \! \dot{\varphi}_t( s_t 
-\pi _t s_t)d\mu _t \, .$$
Thus, up to the action of the non-negative Green's operator, the gradient
of $\kappa$ is given by the right-side of (\ref{evol2}). Indeed, along 
flow paths, the $t$-derivative of $\kappa(\o_t)$ is just the negative
$L_{\o_t}^2$-inner product of $s_t-\pi _t s_t$ and $G_t(s_t-\pi _t s_t)$,
respectively.

\begin{proposition}
Let $\o_t$ be a solution of the initial value problem 
{\rm (\ref{evol})}. Then
$$\frac{d}{dt} \kappa(\o_t) = -\int_M (s_t-\pi _t s_t)
G_t(s_t-\pi _t s_t)
 d\mu_t \, .$$
\end{proposition}

We use this result and the non-negativity of the Green's operator to
rule out non-trivial periodic orbits of the flow (\ref{evol}). This must be
done because in addition to the ones studied above, the flow equation is also 
invariant under the one-parameter group of homotheties, 
where time scales like the square of the distance. In principle, such an 
invariance could give rise to periodic orbits.

\begin{proposition}
The only periodic orbits of the flow equation {\rm (\ref{evol})} are its
fixed points, that is to say, the extremal metrics {\rm (}if any{\rm )} in 
${\mathfrak M}_{\Omega,G}$.
\end{proposition} 

{\it Proof}. Consider the $K$-energy suitably normalized by a volume 
factor to make it scale invariant. If there is a loop solution $\o_t$ of 
(\ref{evol}) for $t\in [t_1,t_2]$, since the volume remains constant, we will
have that $\kappa(\o_{t_1})=\kappa(\o_{t_2})$. By the previous proposition,
since $G_t$ is a non-negative operator, we conclude that 
$G_t(s_t-\pi_t s_t)=0$ on this time interval. This says that $\o_t$ is 
extremal for each $t$ on the interval, and so the 
right side of the evolution equation is zero. Thus, the loop is trivial, a 
fixed point of the flow. \qed

We end this section by showing that the functional (\ref{ene}) decreases
along the flow (\ref{evol}). This should be clear from the way the equation
was set-up, or at the very least, expected. 

\begin{proposition}
Let $\o_t$ be a path in ${\mathfrak M}_{\Omega,G}$ that solves the flow
equation {\rm (\ref{evol})}. Then
$$\frac{d}{dt}\Phi_{\Omega}(\o_t)= -4\int (s_t-\pi_t s_t ) 
L_{t}G_t(s_t-\pi_t s_t ) d\mu_t \leq 0 \, ,$$
and the equality is achieved if and only if $\o_t$ is extremal. In this
expression, $L_t=(\db \d^{\#})^{*}(\db \d^{\#})$ and $G_t$ is the Green's
operator if $\o_t$. 
\end{proposition}

{\it Proof}. Given any variation of the metric with potential function
$\v$, we know that
$$\frac{d}{dt}\Phi(\o_t)=-4 \int s L_t\dot{\v} d\mu_t \, .$$
But $\dot{\v}=G_t(s_t-\pi_t s_t )$, and since $\pi_t s_t$ is a holomorphy
potential and $L_t$ is self-adjoint, we see that
$$\frac{d}{dt}\Phi(\o_t)=-4\int (s_t-\pi_t s_t ) L_{t}G_t(s_t-\pi_t s_t ) 
d\mu_t \, .$$

Both $L_t$ and $G_t$ are non-negative elliptic operators. Thus, $L_tG_t$ is
elliptic and, furthermore, its spectrum is contained in $[0,\infty)$. For if
$L_tG_t \v = \lambda \v$, we have that $G_tL_tG_t \v=\lambda G_t \v$,
and taking the $L^2$-inner product with $\v$ itself, we conclude that
$\lambda$ must be real and non-negative. Therefore, the the expression above 
for the derivative of $\Phi_{\O}$ along flow paths must be non-positive.
If it reaches the value zero at some
$t$, then we must have that $f_t=G_t(s_t-\pi_t s_t )$ is a holomorphy 
potential and $\Delta_t f=(1-\pi_t)s_t$ is an element of the image of 
$1-\pi_t$. Thus, $f_t$ is $L^2$-orthogonal to $(1-\pi_t)s_t$. An
integration by parts argument yields then that $\nabla _t f_t$ must be zero,
and so the function $f_t$ is a constant, which is necessarily zero. 
Thus, $s_t=\pi_t s_t$ and the metric $\o_t$ is extremal. \qed

It is clear that we could have used the function $\Phi_{\O}$ in the r\^ole
that $\kappa$ played when proving that the flow does not have periodic 
orbits other than its fixed points. In fact, it is better to work with
$\Phi_{\O}$ itself. For we do not know if $\kappa$ is in general bounded 
below on ${\mathfrak M}_{\Omega,G}$, but the energy 
functional $\Phi_{\Omega}$ has that property indeed. If the solution
to the flow equation were to exist for all $t\in [0,\infty)$, the 
monotonicity result above would lead us to expect that, as 
$t\rightarrow \infty$, the sequence 
$\o_t$ should be getting closer and closer to an extremal metric. We shall
make several remarks about this possibility later on, but our 
discussion here of both, $\kappa$ and $\Phi_{\Omega}$, serves to show 
the similarities in their behaviour along solution paths to the the extremal 
flow equation. This is rather natural since they are both functions that
can be used to characterize extremal K\"ahler metrics.

\section{The linearized flow equation}\label{sl}
Consider a family of metrics in ${\mathfrak M}_{\Omega,G}$ of the form 
$\o_t(v)=\o_{\v}+i\ddb \a(t,v)$, with $\a(t,0)=0$.
We set $\be=\be_t={\displaystyle \frac{d\a(t,v)}{dv}\mid_{v=0}}$.
The linearization of (\ref{evol2}) at $\o_{\v}$ in the direction of
$\be$ is given by
$$\d _t \be_t =\frac{d}{dv}(G_{(t,v)}((1-\pi_{(t,v)})s_{(t,v)}))\mid_{v=0}\, .
$$
Of course, the argument of the $v$-differentiation in the right side involves 
quantities associated with the metric $\o_t(v)$.

In the remaining part of this section we use the subscript $\v$, or no 
subscript at all, to denote geometric quantities associated with the metric
$\o_{\v}$. 

We have that
$$\frac{d s_{(t,v)}}{dv}\mid_{v=0}=
-\frac{1}{2}\Delta^2_{\v} \be-2(\rho_{\v},i\ddb \be)_{\v}\, .$$
Since the variation of the Green's operator is $-G_\v(\frac{d}{dv}\Delta_{
(t,v)})G_\v$ (keep in mind that this operator needs to be applied only to 
$s-\pi s$, a function that is orthogonal to the constants), 
using the relation between $\rho_\v$ and $\Pi_\v \rho_\v$,
we obtain that 
$$\d _t \be =-\frac{1}{2}\Delta_\v \be -2G_\v(\Pi_\v \rho_\v, i\ddb \be)_\v-
G_\v\left( \frac{d}{dv}\pi_{(t,v)}s_{(t,v)}\mid_{v=0}\right) \, .
$$
By Lemma \ref{le1}, we may write this as
\begin{equation}
\label{lf}
\d _t \be =-\frac{1}{2}\Delta_\v \be -2G_\v(\Pi_\v \rho_\v, i\ddb \be)_\v-
G_\v (\d_{\v}^{\#}\be , X_{\Omega})_{\v}\, , 
\end{equation}
where $X_{\Omega}=\d_{\v}^{\#}(\pi_{\v}s_{\v})$ is the holomorphic vector 
field of the class $\Omega$. Notice that 
$$P_{\v}(\be):=G_\v (\d_{\v}^{\#}\be , X_{\Omega})_{\v}$$ 
is a pseudo-differential operators of order $-1$ in $\be$ whose
coefficients depend non-linearly on the coefficients of the metric $\o_\v$.

We summarize our discussion into the following 

\begin{theorem}
Let $(M,J,\Omega)$ be a polarized K\"ahler manifold and let $G$ be a maximal 
compact subgroup of $Aut(M,J)$. The extremal flow equation 
{\rm (\ref{evol2})} 
{\rm (}or equivalently, {\rm (\ref{evol}))} in 
${\mathfrak M}_{\Omega,G}$ is a non-linear
pseudo-differential parabolic equation.
\end{theorem}

\begin{remark}
For a generic manifold $(M,J)$ the non-trivial holomorphic fields, if any, 
have no zeroes, and the space of real holomorphy potentials reduces to the 
constant functions. Under that hypothesis, the 
pseudo-differential term of order $-1$ in the right side of the linearized 
flow equation (\ref{lf}) vanishes, and the equation reduces to 
$$\d _t \be_t =
-\frac{1}{2}\Delta_\v \be -2G_\v(\Pi_\v \rho_\v,i\ddb \be)_\v \, .$$
This is still a pseudo-differential equation, a zeroth-order 
perturbation of pseudo-differential type of the standard time dependent heat
equation. Thus, even for generic complex manifolds of K\"ahler type, 
the pseudo-differential nature of our flow equation remains in place.
\end{remark}

\begin{remark}
Even when the cohomology class $\Omega$ under study is the canonical 
class $c_1$ (which a fortiori must then have a sign), the extremal 
and Ricci flow do not necessarily coincide with one another.
That will only be the case if we know a priori that $\pi s$ is a constant, 
which as we saw earlier, is a rather non-trivial condition to impose and only 
happens if the Futaki character of the canonical class vanishes. Such a 
restrictive condition fails, for instance, when the manifold is the blow-up of 
${\mathbb C}{\mathbb P}^2$ at one or two points.
\end{remark}
\medskip

We now introduced an {\it approximate} linearized equation whose solution
is needed in our study of local solvability of (\ref{evol2}). In order to
do so, we make some preliminary observations. 

Let $T$ be a positive real number to be determined later and set
$I=[0,T]$. A scale ${\mathcal Y}=\{ {\mathcal Y}_{j}\}_{j\geq 0}$ of 
Banach spaces is a countable family of complete normed spaces such that  
${\mathcal Y}_{j}\supset {\mathcal Y}_{j+1}$ and each ${\mathcal Y}_{j}$
is dense in ${\mathcal Y}_0$. Given one such, we define
$$C_{(j,k)}(I;{\mathcal Y})=C^{0}(I;{\mathcal Y}_{j})\cap \cdots \cap
C^{j-k}(I;{\mathcal Y}_{k}) \, ,$$ 
and provide it with the norm
$$ \| v \| _{j,k}=\sup_{t\in I}\{ \sup_{0\leq r \leq j-k}
\{ \| \partial _{t}^{r}v(t)\|
_{j-r}\} \} \; .$$ 

In what follows, where we shall consider metrics of the form 
$\o_t=\o+i\ddb \v_t$ for path of functions $\v_t$ that begin at $0$ when
$t=0$, we shall always use the scale of Sobolev spaces  
$${\mathcal Y}_{j}=L^{2}_{2j,G}(M)$$
as defined by the background metric $\o$. When $t$ varies on 
the interval $[0,T]$, if we choose $T$ sufficiently small, all the metrics 
$\o_t$ will be equivalent, and the Sobolev spaces defined by them
will be equivalent to each other, with equivalent norms.
We let the Sobolev order jump by $2$ because the operator $F(\v)$ in the 
right side of 
(\ref{evol2}), $$F(\v):=G_{t}(s_t-\pi_t s_t)\, ,$$ 
is of second order, the one reason for the peculiar definition of the scale 
${\mathcal Y}_j$ we shall use. And recall that by the Sobolev embedding 
theorem, we know that $L^{2}_{k}(M)$ is a Banach algebra whenever 
$k>n$. Thus, for as long as the metric $\o_t$ is equivalent to $\o$,
provided that $k>n$, we have a continuous mapping 
$$F: L^{2}_{k+4}(M) \mapsto L^{2}_{k+2}(M)\, .$$

\begin{proposition}\label{cd}
Assume that a solution $\v(t)$ of {\rm (\ref{evol2})} is in 
$C_{(k+1,0)}(I;{\mathcal Y})$ on the interval $I$ for some integer $k$ 
such that $2k>n+2$. Then all the values
of $\d_t^{r}\v(t)$ {\rm (}$1\leq r\leq k+1${\rm )} restricted to $t=0$ 
are completely determined and $\d_t^{r}\v(t)\mid_{t=0}:=\v_r \in {\mathcal
Y}_{k+1-r}=L^{2}_{2k+2-2r,G}(M)$.   
\end{proposition}

{\it Proof}. The initial condition $\v\mid_{t=0}$ is zero, and the equation 
itself sets the value of $\d_t \v \mid_{t=0}=F(0)=
G_{\o}(s_{\o}-\pi_{\o}s_{\o})$ that is evidently in $L^{2}_{2k,G}(M)$.

The relation (\ref{lf}) for $\be=\d_t \v_t$ says that
$$\frac{d}{dt}\be=-\frac{1}{2}\Delta_{\v}\be + P_0(\v)\be \, ,$$
where $P_0$ is a pseudo-differential operator of order zero whose 
coefficients
depend on the coefficients of the metric $\o_\v$ and its curvature
tensor. Since 
$\v_t\in C(I;L^{2}{2k}(M))$ and $2k> n+2$, by the Sobolev embedding theorem,
these coefficients are continuous functions. By regularity of 
pseudo-differential operators on Sobolev spaces, we obtain that
$\d_t \be =\d_t^2 \v_t \in L^{2}_{2(k-1)}(M)$, which is still a continuous 
function because $2(k-1)>n$.

If we differentiate the expression above for $\be=\d_t \v_t$ with respect to 
$t$, we obtain 
$$\frac{d^2}{dt^2}\be=-\frac{1}{2}\Delta_{\v}\frac{d}{dt}\be 
 -\frac{1}{2}L_{\v}(\be)\be+P_0(\v)\frac{d}{dt}\be + P_{0,\v}(\be)\be \, ,$$ 
where $L_{\v}(\be)$ and $P_{0,\v}(\be)$ are the linearizations of 
$\Delta_{\v}$ and $P_{0}(\v)$ at $\v$ in the direction of
$\d_t \v$, respectively. The first is an operator of
order two whose coefficients are continuous. By the metric dependence of
$P_{0}(\v)$, the latter is a pseudo-differential operator of order zero
whose coefficients are also continuous functions. Hence, 
$\d_t^2\be =\d_t^3 \v_t\in L^{2}_{2(k-2)}(M)$. 

Iteration of the argument above yields that
$$\d_t^r \v=F_r(\v, \d_t \v, \ldots, \d_t^{r-1}\v) \, ,$$
where $F_r$ is some operator whose coefficients depend upon the 
coefficients of the metric $\o_\v=\o+i\ddb \v$. 
The desired result for the regularity of $\d_t^r \v$ follows again using 
the Sobolev embedding theorem and the known regularity of the lower
order time derivatives $\d_t^j \v$, $0\leq j\leq r-1$. \qed

Assume given Cauchy data $\v_{0}=0$ for $(\ref{evol2})$
and let $\v_{r}=\partial _{t}^{r}\v(t)\mid _{t=0}$ be the sequence of 
coefficients of the Taylor series  of $\v(t)$ given by the 
proposition above. The Cauchy data $\v_{0}$ determines the sequence 
$\v_{r}$, $1\leq r \leq k+1$. We consider the metric space:
\begin{equation}
W(I)=W^k(I)=\{ \psi(t)\in C_{(k+1,0)}(I;{\mathcal Y}): \; 
\partial _{t}^{r}\psi(t)
\mid _{t=0}=\v_{r}, \; 0\leq r\leq k+1\} \; . 
\label{spa}
\end{equation}
It is not empty, as can be seen by solving the
Cauchy problem for a suitable parabolic equation. 

By a continuity argument, for any $\psi(t)\in W(I)$ the form
$\o_{\psi}=\o + i\ddb \psi(t)$ is positive provided that $t$ is 
sufficiently small. Hence, $\o_{\psi}$ defines a K\"ahler metric.
This metric is not smooth in general. However, if $2k>n+2$, by the Sobolev 
embedding theorem, $\o_{\psi}$ is at least $C^2$, and the operator in the 
right side of (\ref{lf}) will make sense when $\psi$ plays the r\^ole of 
$\v$. Thus, we set
\begin{equation}
P_0(\psi)b=-2G_\psi(\Pi_\psi \rho_\psi, 
i\ddb b)_\psi-G_\psi (\d_{\psi}^{\#}b,X_{\Omega})_{\psi} \, .
\label{lo}
\end{equation}
Then $P_0(\psi)$ is a pseudo-differential operator of order zero in 
$b$, whose coefficients depend upon the coefficients of the 
metric $\o_{\psi}$ and its curvature tensor, all of which are continuous 
functions. For each $t$ on 
a time interval where all the metrics $\o_{\psi}$ are uniformly equivalent, 
we have that 
\begin{equation}
-\frac{1}{2}\Delta_\psi +P_0(\psi) :L^{2}_{2,G}(M) \rightarrow L^{0}_{0,G}(M)=
L^2(M)
\label{ao}
\end{equation}
continuously. We consider the equation
\begin{equation}
\frac{d}{dt}b= -\frac{1}{2}\Delta_{\psi}b +P_0(\psi)b \, ,
\label{ae} 
\end{equation}
whose Cauchy problem will be studied in the next
section. We shall refer to it as the {\it approximate} linearized equation,
the reasons being ---we hope--- clear at this point.

We end this section with the following

\begin{proposition}
Let $\v_1$ be the Cauchy data for {\rm (\ref{ae})}. If
$b(t) \in C_{(k,0)}(I,{\mathcal Y})$ is a solution, then
$\d_t^r b \mid_{t=0}=\v_{r+1}$, $0\leq r \leq k$.
\end{proposition}

{\it Proof}. We have seen above that if $\v(t)$ satisfies (\ref{evol2}), then
$$\d_t^r\v=F_r(\v, \d_t \v, \ldots, \d_t^{r-1}\v) \, , \quad r\geq 2\, ,$$
where $F_r$ is some operator whose coefficients depend upon the 
coefficients of the metric $\o_\v=\o+i\ddb \v$, 
and whose restriction to $t=0$ depends only on the 
sequence $\v_0 , \v_1, \ldots, \v_{r-1}$. The approximate linearized equation
(\ref{ae}) is obtained from the linearization of (\ref{evol2}) given in
(\ref{lf}), when we replace the r\^ole played by $\v(t)$ by that of
$\psi(t)$. But $\psi(t)$ and $\v(t)$ have the same coefficients in their
Taylor expansions up to order $k+1$. Therefore, the solution $b(t)$ to the 
Cauchy problem of (\ref{ae}) with data $b(0)=\v_1$ will have necessarily a 
Taylor series of order $k$ that agrees with the Taylor series of the 
solution to the Cauchy problem of linearized equation (\ref{lf}). The 
conclusion follows by Proposition \ref{cd}. \qed

\section{Local solvability of the extremal flow equation}\label{loc}
In this section, we prove local time existence of solutions to 
the extremal flow (\ref{evol2}). We do so by adapting to our situation
a method of T. Kato for the solvability of abstract differential equations
and non-linear problems \cite{ka}. The pseudo-differential nature of our
linearized equation (\ref{lf}) makes the task harder. But fortunately 
enough, the strictly pseudo-differential part of the equation is lower order,
and most of the analysis is based on that of the standard 
time-dependent heat equation.

\subsection{The Cauchy problem for the approximate linearized equation}
From now on, we take $k$ to be an integer such that $2k> n+2$ and
${\mathcal Y}_{j}=L^{2}_{2j,G}(M)$ as in the previous section.  Given 
Cauchy data $\v_0=0$ for (\ref{evol2}), Proposition \ref{cd} determines
the sequence $\{ \v_{j}\}_{j=0}^{k+1}$, and that in turn allows us to define 
the space $W(I)$ of (\ref{spa}). The interval $I=[0,T]$ will be determined 
later.

For $\psi \in W(I)$, we consider the metrics $\o_{\psi}=\o+i\ddb \psi$ and
the Cauchy problem of the approximate linearized equation (\ref{ae}). Notice
that in that equation, $P_0(\psi)$ is given by (\ref{lo}), a 
pseudo-differential operator of order zero whose coefficients depend 
non-linearly on the coefficients of the metric $\o_{\psi}$ and its curvature 
tensor. 

Let $p(t,s)$ be the evolution operator of 
$$\frac{d}{dt}b=-\frac{1}{2}\Delta_{\psi} b \, .$$
Thus, $p(t,s)$ is a two-parameter family of strongly continuous operators
on ${\mathcal Y}_0$ and ${\mathcal Y}_1$, respectively, such that
$p(t,s)p(s,r) = p(t,r)$, $p(t,t)=1$, and for $b\in {\mathcal Y}_1$ we have
\begin{equation}
\begin{array}{rcl}
\d_t p(t,s)b & = & -\frac{1}{2}\Delta_{\psi(t)}p(t,s)b \, , \\
\d_s p(t,s)b & = & -\frac{1}{2}p(t,s)\Delta_{\psi(s)}b \, .
\end{array}
\label{eop}
\end{equation}
This family of operators exists for $0\leq s\leq t \leq T$, and their 
operator norm is bounded uniformly by a 
constant that only depends upon a bound on $I=[0,T]$ of the coefficients
of $\o_{\psi(t)}$. The function solving (\ref{ae}) with Cauchy data
$\be$ must satisfy the integral equation
\begin{equation}
b(t)=p(t,0)\be+\int_{0}^{t}p(t,s)P_{0}(\psi(s))b(s) d s \, .
\label{ie}
\end{equation}

Consider the set of functions $b(t)$ in $C_{(1,0)}(I;{\mathcal Y})=
C(I;{\mathcal Y}_0)\cap C^{1}(I;{\mathcal Y}_1)$ such that $b(0)=\beta$. The
right hand side of the expression above defines an operator in
this space,
$$P: b \mapsto p(t,0)\beta+\int_{0}^{t}p(t,s)P_{0}
(\psi(s))b(s) d s \, ,$$
and by the explicit form of the coefficients of $P_{0}(\psi)$ mentioned 
above, combined with the continuity of 
pseudo-differential operators on Sobolev spaces, we have that
 $$\| Pb -P\tilde{b}\| \leq C T \| b-\tilde{b}\| \, ,$$
where $C$ is a constant that depends upon the $L^{\infty}$-norm of the 
coefficients of $\o_{\psi(t)}$ and its curvature tensor on the time interval 
$I$. A fixed point argument now yields the following result:

\begin{theorem}\label{thae}
Consider the Cauchy problem for {\rm (\ref{ae})} with Cauchy data 
$b(t)\mid_{t=0}
\in {\mathcal Y}_1$. Then there exists $T$ such that this problem has a
unique solution in $C_{(1,0)}(I;{\mathcal Y})=C(I;{\mathcal Y}_0)\cap 
C^{1}(I;{\mathcal Y}_1)$. The value of $T$ only depends on
supremum norms of the coefficients of $\o_{\psi(t)}$ and its curvature
tensor. 
\end{theorem}

Of course, the regularity of the solution in the theorem above can be
improved if we start with a better initial condition. For that observe 
that the coefficients of the operator $\Delta_{\psi}$ are curves
in $C_{(k,0)}(I;{\mathcal Y})$, and consequently,
$$\Delta_{\psi(t)}: L^{2}_{2j}(M) \mapsto L^{2}_{2j-2}(M)\, , \; 1\leq 
j\leq k\, ,$$
continuously. While the metrics remain equivalent, we can choose a uniform
constant for the operator norm of these maps, and (\ref{eop}) holds for
$b\in L^{2}_{2j}(M)$ with $j$'s in this range. Then we have

\begin{corollary}
If the initial data $b(t)\mid_{t=0}=\v_1 \in {\mathcal Y}_{k}$, the 
solution to the Cauchy problem for {\rm (\ref{ae})} belongs to 
$C_{(k,0)}(I;{\mathcal Y})=C(I;{\mathcal Y}_k)\cap \cdots \cap C^{k}(I,
{\mathcal Y}_0)$.
\end{corollary}
 
{\it Proof}. The arguments in the proof of the theorem and the remarks made 
above show that we now have a solution $b(t)$ to the Cauchy problem for
(\ref{ae}) that is in $C(I;{\mathcal Y}_k)\cap C^{1}(I;{\mathcal Y}_{k-1})$.
This solution satisfies (\ref{ie}) with $\be=\v_1$.

We can differentiate repeatedly the identity (\ref{ae}) in order to show
that the regularity of $b(t)$ with this initial condition can be improved.
Notice that the coefficients of the second order operators 
$d_t^r \Delta_{\psi(t)}$, 
$1\leq r \leq k-1$, are curves in $C(I;L^{2}_{2k-2r})$, and so we have
$d_t^r \Delta_{\psi(t)} \in C(I; {\mathcal L}({\mathcal Y}_{j+r+1}, 
{\mathcal Y}_{j}))$ for $0\leq j \leq k-1-r$. Here, ${\mathcal L}(X,Y)$ is the
space of linear bounded operators from $X$ to $Y$, and the assertion follows
because in the stated range, $L^{2}_{2k-2r}\cdot L^{2}_{2j+2r}\subset 
L^{2}_{2j}$.
This suffices to conclude that the contributions to $d_t^{l+1}b$ arising
from $d_t^{l}\Delta_{\psi}b$ are in $L^{2}_{2k-2l-2}$ if we already 
know that $b\in C_{(k,l)}(I;{\mathcal Y})$.

The analysis of the contributions to $d_t^{l+1}b$ arising from 
$d_t^l (P_0(\psi)b)$ is similar. This time, the coefficients of the
operators $d_t^r(P_{\psi(t)})$ are curves in $C(I;H^{2k-2r-2})$, one 
degree worse than those of $d_t^r \Delta_{\psi(t)}$, but the operators
are of pseudo-differentials of order zero instead. The desired improved 
regularity follows by the same arguments as the ones in the previous 
paragraph. \qed

\subsection{An elliptic equation for $\gamma -F$}
Let us recall that $F(\v)=G_{\v}(s_{\v}-\pi_{\v}s_{\v})$ 
is the second order non-linear operator defined by the right side 
(\ref{evol2}). The derivative $L_{\psi}$ of this map at a general point 
$\psi$ in
${\mathcal Y}_{k+1}$ was computed in \S\ref{sl} and equals the
operator in the right side of (\ref{lf}):
\begin{equation}
L_{\psi}b= -\frac{1}{2}\Delta_\psi b -2G_\psi(\Pi_\psi \rho_\psi, 
i\ddb b)_\psi-G_\psi (\d_{\psi}^{\#}b,X_{\Omega})_{\psi}\, .
\label{lp}
\end{equation}
Since the top part of this
linearization is the negative operator $-\frac{1}{2}\Delta_{\psi}$,
while the lower order term is a pseudo-differential operator of order
zero, coercive estimates for this linearization imply that $\lambda -
L_{\psi}$ is an invertible operator as a map, say, from
${\mathcal Y}_{1}$ to ${\mathcal Y}_{0}$, for a sufficiently large 
constant $\lambda$. 

Let us then take a constant $\lambda$, and consider the non-linear
elliptic map
\begin{equation}
\begin{array}{rcl}
{\mathcal Y}_{k+1} & \longrightarrow & {\mathcal Y}_{k}\\
\v & \mapsto & \lambda \v - F(\v) \, .
\end{array}
\label{ne}
\end{equation}
We remind the reader here of the sequence $\{ \v_{r}\}$ given by 
Proposition \ref{cd}, whose first element is $\v_0=0$. 

\begin{proposition}
For $\lambda$ sufficiently large, there are neighborhoods ${\mathcal O}$ and
${\mathcal V}$ of $\v_0$ and $-\v_1$ in ${\mathcal Y}_{k+1}$ and 
${\mathcal Y}_{k}$, 
respectively, such that the restriction of {\rm (\ref{ne})} to ${\mathcal O}$
is an isomorphism onto ${\mathcal V}$.
\end{proposition}

{\it Proof}. This is a consequence of the Inverse Function Theorem. Indeed,
the linearization $\lambda -L_{0}$ is an invertible operator from
${\mathcal Y}_{1}$ to ${\mathcal Y}_{0}$. Hence, if $f\in {\mathcal Y}_{k}$,
there exists an element $b\in {\mathcal Y}_1$ that satisfies the equation
$$(\lambda-L_{0})b= f \, .$$
We just need now to show that the regularity of $b$ can be improved. 

For that observe that identity above says that the image of $b$ under 
$L_{0}$ is in ${\mathcal Y}_{1}$, and by the regularity properties of 
$\lambda -L_{0}$, we must have
$b\in {\mathcal Y}_2$. By iteration of this argument, we conclude that
$b\in {\mathcal Y}_{k+1}$, and so, $b$ is an element of the tangent space
of ${\mathcal Y}_{k+1}$ at $0$. The desired result follows. \qed 

\begin{corollary}\label{coee}
Let $\psi \in {\mathcal Y}_k$ be sufficiently closed to $-\v_1$. Then, for
large $\lambda$, the equation 
$$\lambda \v - F(\v)= \psi $$
has a solution $\v \in {\mathcal Y}_{k+1}$. The solution is unique if
it is required to be closed enough to $\v_0=0$.
\end{corollary}

In the sequel, we let $D=D^{k}$ be the open neighborhood of $\v_0$ in
${\mathcal Y}_{k+1}$ where the operator $F(\v)$ is defined and smooth.

\subsection{A fixed point argument: local solvability of the non-linear 
equation}
Proceeding by analogy with \cite{ka}, we define $E_{\v_0}(I)$ to be the 
set of curves $\psi(t)\in W^k(I)\subset C_{(k+1,0)}(I;{\mathcal Y})$ such that
$$\| \partial _{t}^{l}\psi(t)-\v_{l}\|_{k+1-l}\leq R,\; l=0,\ldots , k,
\; t\in I \, , $$ for some positive constant $R$. The value of $R$ is chosen
so the ball in ${\mathcal Y}_{k+1}$ with center $\v_0$ and radius
$R$ is contained in the domain $D$ where the operator $F(\v)$ is defined. 
This space is not empty for some $R>0$ and some $I=[0,T]$. 

By the form (\ref{lp}) of the linearization of $F(\v)$ at $\psi$, 
we may conclude that if $\psi_{1}$ and $\psi_{2}$ are elements of 
${\mathcal Y}_{k+1}$, then the operator norm, as a map from ${\mathcal Y}_{k}$
to ${\mathcal Y}_{0}$, satisfies the estimate
$$\| L_{\psi_1}-L_{\psi_2}\| _{k,0} \leq C\| \psi_{1}-\psi_{2}\| _{1}
\, ,$$ for some constant $C$. Indeed, the top part of $L_{\psi}$ in
(\ref{lp}) is half of the 
Laplacian, and its lower order part is a zeroth order pseudo-differential 
operator with nicely behaved coefficients. Then the regularity of 
pseudo-differential operators on Sobolev spaces yields the assertion made.

We now define a key mapping in our proof of the local time
existence to the extremal flow. Let $\psi(t)$ be an element of 
$E_{\v_0}(I)$, and consider the solution $b(t)$ of $(\ref{ae})$ given in 
Theorem \ref{thae}, with initial data $\v_1$. We then solve the 
equation
\begin{equation}
\lambda \v - F(\v) = -b(t)+\lambda \left( \v_0+\int_{0}^{t} b(u) du\right) 
\, , 
\label{abue}
\end{equation}
where we use a real number $\lambda$ such that, if 
$L_0$ is the linearization (\ref{lp}) of $F(\v)$ at $\v=\v_0$, then
$\lambda - L_{0}$ is an isomorphism.

We think of this as a stationary equation in $\v$, that is solved for each 
$t\in I$. Since for $t$ sufficiently small the right side of the equation 
lies in a neighborhood of $-\v_1$, Corollary \ref{coee} applies to produce a 
solution $\v(t)$ in a neighborhood of $\v_0$.

The following two results are the versions of Proposition 7.4 and 
Proposition 7.6 in \cite{ka} adapted to our problem. We give 
proofs here for the sake of completeness.

\begin{proposition}
For sufficiently small $t$, {\rm (\ref{abue})} has a unique solution
$\v(t)$ in a neighborhood of $\v_0$ in  $D \subset {\mathcal Y}_{k+1}$,
with $\v(0)=\v_0 =0$. Furthermore, $\v(t)\in C_{(k+1,0)}(I;{\mathcal Y})$
and $\d_t^r \v(t) \mid_{t=0}=\v_r$, $0\leq r \leq k$ provided $T$ is chosen
sufficiently small, uniformly in $\psi \in E_{\v_0}(I)$. In that case,
$\v(t)\in E_{\v_0}(I)$.
\end{proposition}

{\it Proof}. The operator $\v \mapsto \lambda \v-F(\v)$ is a local 
diffeomorphism of
a neighborhood of $\psi (t)$ in ${\mathcal Y}_{k+1}$ into a 
neighborhood of $\lambda \psi -F(\psi)$ in ${\mathcal Y}_{k}$.
By Theorem \ref{thae}, the right side of (\ref{abue}) is a curve in
$C(I,{\mathcal Y}_k)$ that has value $-\v_1$ at $t=0$. By Corollary
\ref{coee}, we may solve the equation uniquely for $\v(t)$ in a 
${\mathcal Y}_{k+1}$-neighborhood of $\v_0$ and obtain that 
$\v(t)\in C(I,{\mathcal Y}_{k+1})$. This requires to choose $T$
sufficiently small but uniformly in $\psi \in E_{\v_0}(I)$. 

Formal differentiation of the equation solved by $\v(t)$ yields that
$$(\lambda -L_{\v(t)}) \d_t \v =\lambda b - \d_t b=
(\lambda -L_{\psi(t)})b(t) \, .$$
By the invertibility of the operator $\lambda -L_{\v(t)}$ and the known
regularity of the right side, it follows that
$\d_t \v \in C(I,{\mathcal Y}_{k})$ and has value $\v_1$ at $t=0$.
Iterated differentiation yields that $\v(t)\in C_{(k+1,0)}(I;{\mathcal Y})$
and has the desired coefficients in its Taylor series expansion up to 
order $k$. Moreover, the way the equation is solved, we have that
$$\| \d_t^l \v(t) -\v_l\|_{k+1-l} \leq R$$
for $t\in I$. This completes the proof. \qed
 
\begin{proposition}
For $\psi \in E_{\v_0}(I)$, let $\v(t)\in E_{\v_0}(I)$ be the solution 
curve given by the previous proposition. If $T$ is sufficiently small,
the mapping 
\begin{equation}
\begin{array}{c}
E_{\v_0}(I) \longrightarrow E_{\v_0}(I) \\
\psi(t) \longrightarrow \v(t) 
\end{array} \label{fawn} 
\end{equation}
is a contraction in the metric induced by the norm
$||| w |||_{1}=\sup_{t\in I}\| w(t)\| _{1}$, relative to
which, $E_{\v_0}(I)$ is complete.
\end{proposition}

{\it Proof}. Given a curve $b(t)$ in ${\mathcal Y}_q$, we define a norm
by $||| b|||_q=\sup_{t\in I}\| b(t)\|_{q}$. We shall only make use
of the $1$ and $0$ norm, respectively.

Let $\psi_1$ and $\psi_2$ be two elements of
$E_{\v_0}(I)$ and let $b_1$ and $b_2$ be the solutions to the 
corresponding approximate linearized equations with the same initial
condition $\v_1$. We then have that
$$b_1(t)=p_{\psi_1}(t,0)\v_1 \, , \quad b_2(t)=p_{\psi_2}(t,0)\v_1 \, ,$$
where $p_{\psi_1}(t,s)$ and $p_{\psi_2}(t,0)$ are the evolution operators
of the linear equations $\d_t v=L_{\psi_1(t)}v$ and 
 $\d_t v=L_{\psi_2(t)}v$, respectively. Consequently, 
$$b_2(t)-b_1(t)= (p_{\psi_2}(t,0)-p_{\psi_1}(t,0))\v_1 \, ,$$
and using the identity
$$p_{\psi_2}(t,0) \v -p_{\psi_1}(t,0)\v =-\int_{0}^t 
p_{\psi_2}(t,\tau)(L_{\psi_2(\tau)}-L_{\psi_1(\tau)}) p_{\psi_1}(\tau,0)
\v d\tau \, ,$$
we obtain the estimate
$$\| b_2 (t)-b_1(t)\|_{0} \leq C \| \v_1 \|_{k}\int_0^t \| L_{\psi_2(\tau)}-
L_{\psi_1(\tau)}\|_{k,0} d\tau$$
for some constant $C$. But we have observed that
$\| L_{\psi_2(\tau)}-L_{\psi_1(\tau)}\|_{k,0}$ is bounded by a constant
times $\| \psi_2(\tau)-\psi_1(\tau ) \|_{1}$. For small enough $R$, this 
last constant can be chosen uniformly. 
We then obtain that
$$||| b_2 -b_1|||_{0} \leq C T\| \v_1\|_{k} ||| \psi_2 -\psi_1 |||_{1}\, ,$$
showing that the map
$$\psi(t) \mapsto b(t) $$
is a contraction from the $1$-norm to the $0$-norm, with contraction factor 
arbitrarily small with $T$.

That the map $\psi(t) \mapsto \v(t)$ is a contraction now follows because
the map $b(t) \mapsto \v(t)$ is uniformly $C^1$ from the 
$0$-norm to the $1$-norm. This last map is simply the inverse of
$\v \mapsto \lambda \v - F(\v)$ from ${\mathcal Y}_{1}$ to 
${\mathcal Y}_{0}$, and we have that $\lambda - L_{\v(t)}$ is an 
isomorphism from ${\mathcal Y}_{1}$ to 
${\mathcal Y}_{0}$, uniformly in $\psi(t)$ when $\psi(t)$ is close
to $\v_0$. \qed

In view of the previous results, there exists a unique fixed point $\v(t)$
of the map (\ref{fawn}). Since $b(t)$ solves
(\ref{ae}) with initial data $\v_1$, differentiating with respect to $t$ 
in (\ref{abue}) we obtain:
$$(\lambda -L_{\v(t)})\d_{t}\v(t) = -\dot{b}(t) +\lambda b(t)=
(\lambda - L_{\v(t)})b (t)\, ,$$ and since 
$\lambda -L_{\v(t)}$ is injective, we must have that
$$b(t)=\partial _{t}\v(t)\, .$$
We may now use this fact in carrying the time integral in (\ref{abue}), and 
conclude that
$$\frac{d}{dt}\v(t)=F(\v(t)) \, .$$
Thus, the fixed point $\v(t)\in E_{\v_0}(I)$ is a solution to the 
initial value problem (\ref{evol2}). 

We thus arrive at the following 

\begin{theorem}
Let $(M,J,\Omega)$ be a polarized K\"ahler manifold and let $G$ be a maximal 
compact
subgroup of $Aut(M,J)$. The extremal flow equation 
$$\d _t \o _t  = -\rho_t+\Pi _t\rho_t $$
in ${\mathfrak M}_{\Omega,G}$ with a given initial data has a 
unique solution for a short time.
\end{theorem}

In fact, our proof carefully analyses how the time of existence depends 
upon the coefficients of the metric and its curvature tensor. 
Indeed, it shows that the local time of existence depends on
the $L^{\infty}$-norm of the coefficients of the initial metric and
its curvature operator. We can improve a bit the statement above
in relation to the lifespan of the extremal flow. 

\begin{corollary}
Given an initial condition $\o\in {\mathfrak M}_{\Omega,G}$, the extremal 
evolution equation has a unique solution on a maximal time interval 
$0\leq t < T \leq \infty$. If $T< \infty$, then the maximum of the 
point-wise norm of the curvature tensor blows-up as $t\rightarrow T$.
\end{corollary}

The blow-up above, if any, occurs on the Ricci part of the curvature tensor,
rather than the full curvature tensor itself.

\section{Further remarks}
It is of course important to know when the extremal flow has solutions for
all time. Indeed, once the local time existence is known, the next problem to 
consider is the use of the flow to show the existence of extremal 
metrics representing a given cohomology class $\Omega$, task that could be 
accomplished if we manage to prove global time existence and convergence of 
the metrics as $t\rightarrow \infty$. 

This scheme could not possible work in all cases, as we
already know of examples of polarized K\"ahler manifolds without extremal
metrics \cite{bb}. But as a testing ground of its usefulness, we have started 
its analysis when pursuing extremal metrics on polarized manifolds $(M,J,
\Omega)$ with $c_1<0$, or on polarized complex surfaces with $c_1>0$. The 
partial results obtained so far are quite encouraging.

We have two types of fairly strong reasons supporting our belief
that this approach will produce extremal metrics in the said cases.
The first of these is directly related to the flow itself, 
while the other one involves some relation between this flow and the study 
of {\it families} of extremal problems as we vary the cohomology class 
$\Omega$. We discuss them briefly in this section.
 
The evolution equation (\ref{evol}) implies evolution equations 
for various metric tensors associated to the varying metrics. For instance,
the Ricci form evolves according to the equation
$$\frac{d}{dt}\rho=-\frac{1}{2}\Delta \rho + \frac{i\ddb (\pi s)}{2} \, ,$$
the scalar curvature evolves according to the equation
$$\frac{d}{dt}s=-\frac{1}{2}\Delta (s-\pi s) -2 (\rho, i\ddb G(s-\pi s)) \, ,$$
and the Ricci potential evolves according to the equation
$$\frac{d}{dt}\psi =-\frac{1}{2}\Delta \psi -
2G(\rho_H,i\ddb (\psi + G(\pi s))) -(\pi s - s_0) +
\frac{1}{2\mu(M)}\int \psi (s-\pi s)d\mu \, .$$
Here $\rho_H$ is the harmonic component of $\rho$, and $\mu(M)$ is the volume
of $M$ relative to $\o$. 

The first of these equations above shows that the form $\rho$ is a solution 
to the heat equation for the time dependent Hodge Laplacian. One might 
expect that Hamilton's maximum principle (Theorem 9.1 in \cite{ha}) for 
solutions to the heat equation of the rough Laplacian could be extended
to this new setting. If so, such a result would allow us to conclude that
if the initial condition for $\rho$ has a sign, then that sign should
be preserve along the flow (\ref{evol}) for $0\leq t \leq T$, $T$ the 
lifespan of the solution. At the very least, such a result should hold
for generic manifolds $(M,J)$.

We could then apply this to manifolds with no non-trivial holomorphic
vector fields, such as any complex manifold $(M,J)$ with negative first 
Chern class, or most complex surfaces with positive first Chern classes.
Notice that for the blow-up of ${\mathbb C}{\mathbb P}^2$ at one point,
a manifold that carries non-trivial holomorphic vector fields,
the positivity of the Ricci form is preserved along the flow. This makes it 
even more likely that such a result would also hold on any complex surface
with positive $c_1$. 

We may also refine our earlier Theorem \ref{tc} when dealing with a 
complex surface of positive first Chern class. Indeed, we have the following
result, whose proof will be given elsewhere.

\begin{theorem}
Let $(M,J,\O)$ be a polarized complex surface of positive first Chern class.
Given any K\"ahler metric $g$ in ${\mathfrak M}_{\O,G}$, the image of the
holomorphy potential $\pi_g s_g$ is an interval contained in the set of 
positive real numbers, interval that only depends on $\O$ and not on $g$.
\end{theorem}

Thus, if for a given initial condition with positive Ricci
curvature on a $c_1$-positive surface we have that solutions to the flow 
(\ref{evol}) exists for all time, and converge to an extremal metric as
time goes to infinity, the extremal metric so obtained
would have positive scalar curvature, as expected.

The preservation of the sign of the Ricci tensor should have very strong 
implications on the global analysis of (\ref{evol}). This  property has 
been of utmost importance already in the work of Hamilton \cite{ha}, and 
should remain so in the general analysis of our flow equation as well. 
If $c_1>0$, we could combine this with plausible global time 
existence results, and pass to a Cheeger-Gromov-Hausdorff 
limit, an important step towards settling the convergence issue. 

We venture the following two conjectures.

\begin{conjecture}
Let $(M,J)$ be a complex manifold of K\"ahler type polarized by a 
K\"ahler class $\Omega$. If $c_1(M,J)<0$, there exists an initial condition
to the extremal flow {\rm (\ref{evol2})} equation so that the solution 
exists on $[0,\infty)$ and, as $t\rightarrow \infty$, converges to a metric 
of constant negative scalar curvature representing $\Omega$. 
\end{conjecture}
\medskip

\begin{conjecture}
Let $(M,J)$ be a complex surface of positive first Chern class polarized by
a K\"ahler class $\Omega$. Then there exists an initial condition
to the extremal flow {\rm (\ref{evol2})} equation so that the solution 
exists on $[0,\infty)$ and, as $t\rightarrow \infty$, converges to an 
extremal metric of positive scalar curvature representing $\Omega$.
\end{conjecture}
\medskip

The initial condition we have in mind in these two cases is given by a
metric whose Ricci form is negative or positive, respectively. After the
work of Yau \cite{ya} on the Calabi conjecture, we know we can always 
find this type of metrics on any given polarization.

These conjectures are further supported by the results in \cite{ss2}, that
we proceed to describe in brief detail. For a complex manifold
$(M,J)$ of complex dimension $n$, we denote by ${\mathfrak M}$ the space 
of K\"ahler metrics on $(M,J)$. As before, given a positive class 
$\Omega \in H^{1,1}(M,{\mathbb C})\cap 
H^{2}(M,{\mathbb R})$, we let ${\mathfrak M}_{\Omega}$ be the space of 
of K\"ahler metrics whose K\"ahler forms represent $\Omega$. We shall
also consider the space ${\mathfrak M}_1$ of 
K\"ahler metrics of volume one, and ${\mathcal K}_{1}$, the space of 
cohomology classes that can be represented by K\"ahler forms of metrics in
${\mathfrak M}_1$:
\begin{equation}
{\mathcal K}_{1}=\{ \Omega \in H^{1,1}(M,{\mathbb C}): \; \Omega= [\o]
\; {\rm for \; some\; }\o \in {\mathfrak M}_{1}\} \, .\label{cvo}
\end{equation}

Extremal metrics in ${\mathfrak M}_{\Omega}$ achieve the infimum of the 
functional $\Phi_{\Omega}$ in (\ref{ene}), and we have the lower bound
(\ref{lb}):
$$E(\Omega)=\int (\pi _g s_g)^2 d\mu_g \, .$$

One approach to providing $(M,J)$ with a canonical shape would be to 
find critical points of the functional 
\begin{equation}
\begin{array}{rcl}
{\mathfrak M}_1 & \rightarrow & {\mathbb R} \\
\omega & \mapsto & {\displaystyle \int_M s_{\omega}^2 d\mu _{\omega }}\, .
\end{array} \label{fun1}
\end{equation}
A special metric $\o$ of this type must have the following properties:
\begin{enumerate}
\item[a)] $\o$ achieves the lower bound $E([\o])$, that is to say, $\o$ is
extremal relative to the polarization defined by the K\"ahler class $\Omega=
[\omega]$ that it represents; 
\item[b)] the K\"ahler class $\Omega=[\omega]$ is a critical point of 
$E(\Omega)$ 
as a functional defined over ${\mathcal K}_1$.
\end{enumerate}
Thus, the search for critical points of (\ref{fun1}) ---or {\it strongly
extremal metrics} \cite{si2}--- achieving an optimal lower bound involves the 
solution of back-to-back minimization problems: the
first solving for critical points of (\ref{ene}) within a fixed cohomology
class $\Omega$, and the second solving for those classes that minimize the
critical value $E(\Omega)$ as the class $\Omega$ varies within 
${\mathcal K}_1$. Naturally, we separate the two problems by, in addition
to (\ref{ene}), introducing the functional

\begin{equation}
\begin{array}{rcl}
{\mathcal K}_1 & \rightarrow & {\mathbb R} \\
\Omega & \mapsto & E(\Omega)=
{\displaystyle \int_M (\pi s)^2 d\mu }\, ,
\end{array} \label{cun1}
\end{equation}
where the geometric quantities in the right are those associated with any
$G$-invariant metric that represents $\Omega$, for $G$ a fixed maximal 
compact subgroup of the automorphism group of $(M,J)$. Its extremal 
points will be called either critical or canonical classes. We then have 
\cite{ss} the following 

\begin{theorem} \label{count}
Let $\Omega$ be a cohomology class that is represented by a K\"ahler
metric $g$, assumed to be invariant under the maximal compact 
subgroup $G$ of the biholomorphism group of $(M,J)$. Then $\Omega$ is 
critical class if and only if
$$\int_{M} (\pi_g s_g) (\Pi_g \rho , \a) d\mu_g =0 $$
for any trace-free harmonic $(1,1)$-form $\a$. In this expression, 
$\rho$ is the Ricci form of the metric $g$, 
$\pi$ is the $L^2$ projection {\rm (\ref{proj1})} onto the space
of holomorphy potentials, and $\Pi$ is its lift {\rm (\ref{proj})} at the 
level of {\rm (1,1)}-forms.
\end{theorem}

This theorem states that $\Omega$ is a critical class of 
(\ref{cun1}) if and only if
$$\int_{M} (\pi_g s_g) (\Pi_g \rho , \a) d\mu_g =0 $$
for any trace-free harmonic $(1,1)$-form $\a$. In other words, 
the form $\pi s \Pi \rho$ is $L^2$-perpendicular to the space of 
trace-free harmonic (1,1)-forms, and therefore, by Hodge decomposition, the
class must be such that
\begin{equation}
\pi s \Pi \rho = \lambda \omega +\partial G_{\partial}(\partial^{*}(\pi s\,
 \Pi \rho)) +\partial^{*}G_{\partial}(\partial(\pi s\, \Pi \rho))\, ,
\end{equation}
for $\lambda$ equal to the $L^2$-projection of $(\pi s)^2$ onto
the constants, divided by $2n$:
\begin{equation}
\lambda = \frac{1}{2n}\int (\pi s)^2 d\mu_g \, . \label{lam}
\end{equation}

In order to study the existence of critical classes, we may 
consider \cite{ss2} the evolution equation 
\begin{equation}
\frac{d\Omega}{dt} = \pi s \, \Pi \rho - \lambda \omega +
\partial G_{\partial}(\partial^{*}(\pi s\,
 \Pi \rho)) +\partial^{*}G_{\partial}(\partial(\pi s\, \Pi \rho))\, .
\label{hfe}
\end{equation}

The flow equation (\ref{hfe}) defines a dynamical system on
${\mathcal K}_1$ provided the solutions remain in ${\mathcal K}_1$ 
throughout time. Unfortunately, this is not true in general \cite{ss2}.

In the generic case where all non-trivial holomorphic vector fields of
$(M,J)$ have no zeroes, equation (\ref{hfe}) can be extended to a dynamical 
system on
$$\overline{{\mathcal K}}_1=\{ \Omega \in H^{1,1}(M,{\mathbb C})\cap
H^2(M,{\mathbb R}): \; \frac{\Omega^n}{n!}=1\} \, .$$
Indeed, given $\Omega \in 
\overline{{\mathcal K}}_1$, let us define the function
$$s_{\Omega} := 4\pi n \frac{c_1 \cdot \O ^{n-1}}{\O^n} 
\, .$$ 
If $\Omega$ were a K\"ahler class represented by a metric $g$, this function 
would be precisely the holomorphy potential $\pi_g s_g$.
The equation
\begin{equation}
\frac{d}{dt}\O= 2\pi s_{\Omega} c_1 - \frac{s_{\Omega}^2 }{2n}\O \, , 
\label{hfe2}
\end{equation}
extends (\ref{hfe}), which as such is defined only on ${\mathcal K}_1$,
all the way to a dynamical system on $\overline{{\mathcal K}}_1$.

Solutions to (\ref{hfe2}) with initial data
in $\overline{{\mathcal K}}_1$ remain in $\overline{{\mathcal K}}_1$. 
In fact, we have that \cite{ss2}

\begin{theorem}
Suppose that all non-trivial holomorphic vector fields of $(M,J)$ 
have no zeroes. 
Then solutions to {\rm (\ref{hfe2})} with initial data in 
$\overline{{\mathcal K}}_1$ converge, as $t\rightarrow \infty$, to a 
stationary point of the equation in the space $\overline{{\mathcal K}}_1$. 
\end{theorem}

It is then of natural interest to see if solutions to the equation with
Cauchy data given by a positive class, that is to say, an element of
${\mathcal K}_1$, remain positive thereafter. We already know \cite{ss2} of 
examples where this is not so, with solutions to the flow
equation that are initially in the K\"ahler cone but that, in converging to a 
critical point of the flow in $\overline{{\mathcal K}}_1$, must 
eventually leave the cone through its walls. 

In fact, this situation occurs already on complex surfaces, where the 
stability of 
${\mathcal K}_1$ under the flow
(\ref{hfe2}) can be analyzed using a criterion giving necessary
and sufficient for a cohomology class to be K\"ahler, criterion that
extends that of Nakai for integral classes.
Applied to our problem, if the Chern number $c_1^2 \neq 0$, we have that a 
path $\O_t$ solving (\ref{hfe2}) with initial condition in
${\mathcal K}_1$ stays there forever after if, and only if,
$$\O_0 \cdot [D] +8\pi^2 (c_1\cdot \O_0)( c_1\cdot[D])
\left( \frac{e^{c_1^2 t}-1}{c_1^2}\right) > 0$$
for all $t\geq 0$ and for all effective divisors $D$ in $(M,J)$.
When $c_1^2=0$ we still obtain a similar criterion, replacing the expression 
in parentheses above by its limit $t$ as $c_1^2 \rightarrow 0$.

This forward stability of the K\"ahler cone holds in very general 
situations, as can be seen by a run-down of the various cases in the 
Enriques-Kodaira classification of complex surfaces \cite{ss2}. In 
particular, it holds 
if the complex surface has a signed first Chern class $c_1$, 
condition under
which all solutions to the flow (\ref{hfe2}) that start in
${\mathcal K}_1$ stay there forever after, and as $t\rightarrow \infty$,
they either converge to the only critical class 
$\sqrt{2}({\rm sgn}\, c_1)c_1/c_1^2$ of (\ref{cun1}) if $c_1^2>0$, or all 
classes are critical and the flow is constant if $c_1=0$.

Notice that the positivity condition above involves the evaluation of $c_1$ 
over the divisor $D$, and only in the case when there are effective 
divisors $D$ for which $c_1\cdot[D]$ changes sign from one to another could 
the condition fail to hold. Merely fixing the sign of $c_1$ prevents this
from happening, but the counterpart to that is of great interest.
It shows that the existence of divisors on which $c_1$ achieves values
of opposite signs is in effect part of the reason why the
the K\"ahler cone might be poorly behaved in relation to 
the flow (\ref{hfe}).

When the surface in question has positive first Chern class and
carries non-trivial holomorphic fields, the forward stability of the 
K\"ahler cone under the flow (\ref{hfe}) seems to hold also, though we have
only verified that for the case of ${\mathbb C}{\mathbb P}^2$ blown-up
at one point.

In higher dimension and for manifolds $(M,J)$ where $c_1$ is either 
positive or negative, the space of K\"ahler classes is also forward stable 
under the flow (\ref{hfe2}). As a matter of fact, there is a positivity 
criterion that generalizes the one outlined above for surfaces, which 
guarantees forward stability of the K\"ahler cone under the flow. 
Manifolds with signed first Chern classes meet this criterion, though 
for these particular cases one can also give a direct argument that proves the 
flow stability of the cone.

All of these facts combined give further support to the conjectures
made earlier. We end up venturing a final one. 

\begin{conjecture}
Suppose the flow equation {\rm (\ref{hfe2})} with initial data in the 
K\"ahler cone converges to a stationary point that is outside it.
Then the extremal K\"ahler cone is not a closed subset of the K\"ahler cone.
\end{conjecture}
\medskip

In other words, under the given hypothesis, there should exist cohomology 
classes in the K\"ahler cone that cannot be represented by extremal metrics.

\end{document}